\documentclass[11pt,a4paper,oneside,reqno]{amsart}
\usepackage{amsfonts}
\usepackage{amssymb}
\usepackage{mathrsfs}
\usepackage[super]{nth}
\usepackage{amsthm}
\usepackage{soul}
\usepackage[defaultcolor=orange]{changes}
\usepackage{amsmath}
\usepackage{enumitem}
\usepackage{mathtools}
\usepackage[bottom]{footmisc}
\usepackage{pifont}
\usepackage{float}
\usepackage{pgfplots}
\usepackage{color}
\usepackage{url}
\usepackage{tikz-cd}
\usepackage{graphicx}
\usepackage{csquotes}
\usepackage[top=1in, bottom=1in]{geometry}
\graphicspath{ {images/} }

\usepackage[hidelinks]{hyperref}
\usepackage{makecell}
\usepackage{fancyhdr}

\newcommand\cyr
{
\renewcommand\rmdefault{wncyr}
\renewcommand\sfdefault{wncyss}
\renewcommand\encodingdefault{OT2}
\normalfont
\selectfont
}
\DeclareTextFontCommand{\textcyr}{\cyr}

\newlist{case}{enumerate}{1}
\setlist[case,1]{
  label={\textsc{Case}~\arabic*:},
  leftmargin=*,
  align=left,
  labelsep=5mm
}

\newlist{step}{enumerate}{1}
\setlist[step,1]{
  label={\textsc{Step}~\arabic*:},
  leftmargin=*,
  align=left,
  labelsep=5mm
}

\allowdisplaybreaks
\usepackage[OT2,OT1]{fontenc}
\usepackage{wasysym}
\DeclareSymbolFont{extraup}{U}{zavm}{m}{n}
\DeclareMathSymbol{\varheart}{\mathalpha}{extraup}{86}
\DeclareMathSymbol{\vardiamond}{\mathalpha}{extraup}{87}
\DeclareMathSymbol{\varspade}{\mathalpha}{extraup}{85}
\pgfdeclarelayer{nodelayer}
\pgfdeclarelayer{edgelayer}
\pgfsetlayers{nodelayer,edgelayer,main}
\def \proof {\noindent\textit{Proof.}\hspace{5mm}}

\def \bs {$\hfill \blacksquare$\\}
\def \bsn {$\hfill \blacksquare$}
\def \sq {$\hfill \square$\\}
\def \sqn {$\hfill \square$}

\DeclareMathOperator{\supp}{supp}

\DeclareMathOperator{\charac}{char}
\DeclareMathOperator{\Z}{\mathbb{Z}}
\DeclareMathOperator{\N}{\mathbb{N}}
\DeclareMathOperator{\Q}{\mathbb{Q}}

\DeclareMathOperator{\R}{\mathbb{R}}
\DeclareMathOperator{\C}{\mathbb{C}}

\def\forkindep{\mathrel{\raise0.2ex\hbox{\ooalign{\hidewidth$\vert$\hidewidth\cr\raise-0.9ex\hbox{$\smile$}}}}}
\definecolor{supergreen}{RGB}{0, 170, 0}
\definecolor{superred}{RGB}{170, 0, 0}
\definecolor{darkred}{RGB}{100, 0, 0}

\newcommand{\Th}{\mbox{\ttfamily Th}}

\newcommand{\TTh}{\mbox{\emph{\ttfamily Th}}}

\newcommand{\defeq}{\mathrel{\mathop:}=}

\newtheorem{thm}{Theorem}[section]
\newtheorem{prop}[thm]{Proposition}
\newtheorem{cor}[thm]{Corollary}

\newtheorem{lem}[thm]{Lemma}

\newtheorem{prob}[thm]{Problem}

\theoremstyle{definition}
\newtheorem{definition}[thm]{Definition}
\newtheorem{exmp}[thm]{Example}
\newtheorem{remark}[thm]{Remark}
\newtheorem{innercustomthm}{Theorem}
\newenvironment{customthm}[1]
  {\renewcommand\theinnercustomthm{#1}\innercustomthm}
  {\endinnercustomthm}

\newtheorem{innercustomcor}{Corollary}
\newenvironment{customcor}[1]
  {\renewcommand\theinnercustomcor{#1}\innercustomcor}
  {\endinnercustomcor}
\newtheorem{innercustomconj}{Conjecture}
\newenvironment{customconj}[1]
  {\renewcommand\theinnercustomconj{#1}\innercustomconj}
  {\endinnercustomconj}

\newtheorem{innercustomass}{Assumption}
\newenvironment{customass}[1]
  {\renewcommand\theinnercustomass{#1}\innercustomass}
  {\endinnercustomass}

\pgfplotsset{compat=1.14}

\begin{document}
\title{Finite Undecidability in Fields I: NIP Fields}
\author{Brian Tyrrell}
\address{Mathematical Institute, Woodstock Road, Oxford OX2 6GG.}
\email{brian.tyrrell@maths.ox.ac.uk}

\begin{abstract}
A field $K$ in a ring language $\mathcal{L}$ is \emph{finitely undecidable} if $\mbox{Cons}(T)$ is undecidable for every nonempty finite $T \subseteq \mbox{\ttfamily Th}(K; \mathcal{L})$. We extend a construction of Ziegler and (among other results) use a first-order classification of Anscombe \& Jahnke to prove every NIP henselian nontrivally valued field is finitely undecidable. We conclude (assuming the \textit{NIP Fields Conjecture}) that every NIP field is finitely undecidable. This work is drawn from the author's PhD thesis \cite[Chapter 3]{tyrrellphd}.
\end{abstract}

\maketitle
\vspace{-2mm}
\section{Introduction}
\setul{1.5pt}{.4pt}

\noindent The author was motivated to consider this topic by the following question:

\begin{prob}\label{1}
Does there exist an infinite, finitely axiomatisable field? 
\end{prob}

This (still open) problem was posed explicitly by I.\ Kaplan at the 2016 Oberwolfach workshop on \emph{Definability and Decidability Problems in Number Theory} \cite[Q4]{ober}, though existed as folklore before. This relates closely to another elementary question:

\begin{prob}\label{2}
Does there exist a finitely axiomatisable theory of fields which is decidable and has an infinite model?
\end{prob}

That is, this problem is to ascertain the existence of an infinite field $F$ and a collection of first-order sentences $T$ in the language of rings $\mathcal{L}_r$, such that $T \subseteq \Th(F; \mathcal{L}_r)$, $T$ is finitely axiomatised, $T$ models the field axioms, and there exists a decision procedure to determine membership of $\mbox{Cons}(T)$ -- the set of $\mathcal{L}_r$-sentences $\phi$ such that $T \models \phi$.

Answering \emph{Problem \ref{2}} in the negative (as the empirical evidence suggests might indeed be the case) would answer \emph{Problem \ref{1}} in the negative too. This is the focus of modern investigations. One approach to this was established by Ziegler \cite{ziegler}, and generalised further by Shlapentokh \&\ Videla \cite{shlapentokhvidela}. Ziegler's idea was to take a finitely axiomatised subtheory of a sufficiently saturated field with a powerful model completeness property (he considered $\C$, $\widetilde{\mathbb{F}_p(t)}$ -- where `$\,\widetilde{F}\,$' denotes the \emph{algebraic} closure of a field $F$ -- $\R$, and $\Q_p$) and prove it to be a subtheory of a field interpreting arithmetic. By using a result of Tarski (see below), there is no nonempty finitely axiomatised subtheory (``finite subtheory'') of {\ttfamily ACF${}_0$}, {\ttfamily ACF${}_p$}, {\ttfamily RCF}, or {\ttfamily $p$CF} that is decidable. With this in mind we forward the following definition\footnote{Shlapentokh \& Videla \cite{shlapentokhvidela} call this property \emph{finite hereditary undecidability}; for notational ease we remove the word ``hereditary''.}:

\begin{definition}
A theory $T$ in a language $\mathcal{L}$ is \emph{finitely undecidable} if every finitely axiomatised $\mathcal{L}$-subtheory of $T$ is undecidable. (An $\mathcal{L}$-structure is \textit{finitely undecidable} if its complete $\mathcal{L}$-theory is.)
\end{definition}

What infinite fields are finitely undecidable (implicitly: in $\mathcal{L}_r$)? That is the motivating question for this paper and its sequel. If there exists a field satisfying \emph{Problem \ref{1}} or \emph{\ref{2}}, this field is not finitely undecidable; we will show however there is a considerably broad class of fields whose members have this undecidability property. 

In particular, most (if not all) infinite fields whose model theory in the language of rings is well understood will be finitely undecidable, as we will argue. Following this philosophy, we are motivated by a long-standing model theory conjecture:

\begin{customconj}{(Shelah's NIP Fields Conjecture)}
\emph{Every infinite NIP field is either separably closed, real closed, or admits a nontrivial henselian valuation.}
\end{customconj}

In \cite{btyrrelpac} we consider other classification-theoretic conjectures on field theories. Here, we adapt Zielger's argument to {\ttfamily SCF}${}_{p, e}$ ($\S \ref{scf}$) and the complete theories of certain henselian valued fields in the language $\mathcal{L}_{val}$; the language of rings with an additional unary predicate for the valuation ring (\S \ref{hvfsect} \& \S\ref{hvfagain}). Our main result is:

\begin{customthm}{(\ref{thisis} + \ref{newmixed} + \ref{nipequichar})}
\textit{Let $(K, v)$ be an} 
\begin{itemize}
    \item \textit{equicharacteristic $0$, or}
    \item \textit{mixed characteristic, or}
    \item \textit{equicharacteristic $p>0$ separably defectless Kaplansky}
\end{itemize}
\textit{henselian nontrivially valued field. Then $\TTh(K; \mathcal{L}_{val})$ is finitely undecidable.}

\textit{Moreover, if $\mathcal{O}_v$ is $\mathcal{L}_r$-definable, $K$ is finitely undecidable} as a field.
\end{customthm}

Thanks to some deeper results in classification theory, some immediate consequences of this are:

\begin{customcor}{\ref{finalcor}}
\phantom{boo}
\begin{enumerate}
    \item[\textit{(1)}] \textit{If $(K, v)$ is an NIP henselian nontrivially valued field, $\TTh(K; \mathcal{L}_{val})$ is finitely undecidable. Furthermore if $\mathcal{O}_v$ is $\mathcal{L}_r$-definable in $K$, then $K$ is finitely undecidable as a field.}
    \item[\textit{(2)}] \textit{Every infinite \emph{dp-finite} field is finitely undecidable.}
    \item[\textit{(3)}] \textit{Assuming the \emph{NIP Fields Conjecture}, {every} infinite NIP field is finitely undecidable.} 
\end{enumerate}
\end{customcor}

We recall the notion of a \emph{dp-finite} field $L$ below (for a full, general discussion of dp-rank, see \cite{kaplan, simon}). We work in a sufficiently saturated model of $\Th(L; \mathcal{L}_r)$.

\begin{definition}
For $n \in \N_{>0}$, let $I_1, \dots, I_n$ be a list of sequences, and $A$ a set of parameters. We say that the sequences $I_1, \dots, I_n$ are \textit{mutually indiscernible over $A$} if for each $1 \leq t \leq n$, the sequence $I_t$ is indiscernible over $A \cup I_{\neq t}$.

If $p$ is a partial type over $A$, we say \textit{$p$ is of dp-rank $\geq n$} if there exists $a \models p$ and a list of sequences $I_1, \dots, I_n$ mutually indiscernible over $A$ such that $I_k$ is not indiscernible over $A \cup \{a\}$ for $1 \leq k \leq n$. (The dp-rank of a partial type is always $\geq 0$.)
\end{definition}

This definition does not in fact depend on the parameter set $A$, as evidenced by Simon \cite[Lemma 4.14]{simon}.

\begin{definition}
A field $L$ is \emph{of finite dp-rank} if for some $m \in \N$, ``$x = x$'' has dp-rank $\geq m$ but does not have dp-rank $\geq m+1$.
\end{definition}

See \emph{Example \ref{finalex}} for `interesting' examples of \emph{Theorem $($\ref{thisis} $+$ \ref{newmixed} $+$ \ref{nipequichar}$)$ \& Corollary \ref{finalcor}} in action. 

All of our undecidability results rely on a theorem of Tarksi, which we now state: general references for this material are \cite{eltt, shoenfield, tmr}. 

\begin{thm}\label{a2}
\emph{\textbf{(Tarski)}} Let $\mathcal{L}_1$, $\mathcal{L}_2$ be finite languages. The $\mathcal{L}_2$-theory $T_2$ is hereditarily undecidable if there exists a finitely axiomatised essentially undecidable $\mathcal{L}_1$-theory $T_1$ and models $M_1 \models T_1$, $M_2 \models T_2$ such that $M_1$ is interpretable in $M_2$.
\end{thm}

\proof
{This is \cite[pp.\ 87--89]{eltt}, using different (but equivalent) terminology. Results such as this originate in \cite[\S I.3 \& \S I.4]{tmr}; cf.\! \textit{Theorems 6, 7 \& 8 ibid.}} \bs

\noindent We will use this as follows:

\begin{lem}\label{prop112}
\emph{\cite[Proposition 11.2]{shoenfield}.} Let $\mathcal{L}$ be a finite language and $a_1, \dots, a_n$ constant symbols not in $\mathcal{L}$. Let $M$ be an $\mathcal{L}(a_1, \dots, a_n)$-structure and $M|_{\mathcal{L}}$ the reduct of $M$ to $\mathcal{L}$. If $\TTh(M; \mathcal{L}(a_1, \dots, a_n))$ is hereditarily undecidable, so too is $\TTh(M|_{\mathcal{L}}; \mathcal{L})$. \bsn
\end{lem}

\begin{cor}\label{appendixacor}
Let $K$ be a field of characteristic $0$, and $L$ a field of characteristic $p > 0$ such that there exists $t \in L$ transcendental over $\mathbb{F}_p$. Let $\mathcal{L}$ be a finite expansion of the language of rings. Suppose $\Z$ is $\mathcal{L}$-definable with parameters in $K$ and $\mathbb{F}_p[t]$ is $\mathcal{L}$-definable with parameters in $L$. Then $\TTh(K;\mathcal{L})$ and $\TTh(L; \mathcal{L})$ are hereditarily undecidable.
\end{cor}

\proof
This is an application of \emph{Theorem \ref{a2}} with $\mathcal{L}_2$ an expansion of $\mathcal{L}$ by constant symbols, $\mathcal{L}_1 = \mathcal{L}_r$, $T_1 = Q$ (\emph{Robinson arithmetic}) and $M_1 = \N$. Notice $\N$ is $\emptyset$-$\mathcal{L}_r$-definable in $\Z$ (e.g.\ by \emph{Lagrange's Four Square Theorem}) hence interpretable in $\Z$. By \cite[\S 4a--4b]{undecrob}, $\N$ is interpretable in the $\mathcal{L}_r(t)$-structure $\mathbb{F}_p[t]$. Thus, by assumption $\N$ is interpretable in $K$ (with, say, parameters $\overline{c} = \{c_1, \dots, c_k\}$) and $\N$ is interpretable in $L$ (with, say, parameters $\overline{d} = \{d_1, \dots, d_l\}$). By \emph{Theorem \ref{a2}}, $\Th(K; \mathcal{L}(\overline{c}))$ and $\Th(L; \mathcal{L}(\overline{d}))$ are hereditarily undecidable. The hereditary undecidability of $\Th(K; \mathcal{L})$ and $\Th(L; \mathcal{L})$ follows from \emph{Lemma \ref{prop112}} exactly. \bs

\begin{remark}
\textit{Notational remark}. By ``$K \equiv_{\mathcal{L}} L$'' we denote that two $\mathcal{L}$-structures $K, L$ are elementarily equivalent. If ``$K \equiv L$'' is written, the language is implicitly $\mathcal{L}_r$. In the case of \emph{valued fields}, we will frequently write ``$(K, v) \equiv (L, w)$'' to denote an $\mathcal{L}_{val}$-elementary equivalence, where $v, w$ are valuations on their respective fields. \sqn
\end{remark}

\bigskip
\section{First Exploration by Ziegler}
Ziegler's main result of \cite{ziegler} is the construction of a field $K_q$ satisfying the following theorem:

\begin{thm}\label{ziegmainthm}
Let $L$ be the field $\C$, $\widetilde{\mathbb{F}_p(t)}$, $\R$, or $\Q_p$, and $q \neq p$ prime. There exists a field $K_q \subseteq L$ such that
\begin{enumerate}
    \item {$\Z$ or $\mathbb{F}_p[t]$ is definable (with parameters)} in $K_q$;
    \item If the intermediate field $K_q \subseteq H \subseteq L$ is finite over $K_q$, either $[H : K_q] = 1$ or $q | [H : K_q]$.
\end{enumerate}
\end{thm}

{\proof
This is \cite[Theorem, p.~270]{ziegler} with a specified $A$. Cf.\ \textit{Theorems 3.3 \& 4.4,} and \S 5, of \cite{shlapentokhvidela} (this paper, by Shlapentokh \& Videla, generalises Ziegler's construction). \bs}

\begin{cor}\label{zieglercor}
\emph{\textbf{(Ziegler)}} Let $p$ be prime. \emph{{\ttfamily ACF}${}_0$, {\ttfamily ACF}${}_p$, {\ttfamily RCF}}, and \emph{{\ttfamily $p$CF}} are finitely undecidable.
\end{cor}

\proof
{\cite[Corollary, p.~270]{ziegler}}. Let $T$ be a finite subtheory of $L$, which is one of the fields $\C$, $\widetilde{\mathbb{F}_p(t)}$, $\R$, $\Q_p$. Let $P$ be the set of primes distinct to $p$. For each $q \in P$, use \emph{Theorem \ref{ziegmainthm}} to obtain a field $K_q$ satisfying \emph{(1)} \&\ \emph{(2)}. Let $\mathcal{U}$ be a nonprincipal ultrafilter on $P$, and let $\mathbb{K}$ be the ultraproduct $\prod_{q \in P} K_q / \mathcal{U}$.

We claim that $\mathbb{K}$ is relatively algebraically closed in $L^{\mathcal{U}}$. Indeed, suppose $f = (f_q) \in \mathbb{K}[X]$ has a root $\alpha = (\alpha_q) \in L^{\mathcal{U}}$. {In which case, $\{q\mbox{ : } f_q(\alpha_q) = 0\} \in \mathcal{U}$, hence $\{q\mbox{ : } [K_q(\alpha_q) : K_q] \leq \deg(f)\} \in \mathcal{U}$. Consequently, as $q | [H : K_q]$ for all proper finite extensions $K_q \subset H \subseteq L$, $\{q\mbox{ : } [K_q(\alpha_q) : K_q] =1 \} \in \mathcal{U}$. Hence $\alpha \in \mathbb{K}$ as desired.}

By the model theory of algebraically/real closed/$p$-adically closed fields\footnote{See {\cite[\S 3.2 \& 3.3]{marker}} as a reference for the former two, and {\cite[\S 5]{prestelandroq}} for the latter.}, we deduce $\mathbb{K} \equiv L^{\mathcal{U}}$ ($\equiv L$ by \textit{\L o\'{s}' Theorem}; specifically \cite[Theorem 4.1.9]{changkeisler}). Therefore $\mathbb{K} \models T$. As $T$ is finitely axiomatised, by \textit{\L o\'{s}' Theorem} there must exist some $q \in P$ such that $K_q \models T$. Thus by \emph{Theorem \ref{ziegmainthm} \& Corollary \ref{appendixacor}}, $\Th(K_q; \mathcal{L}_r)$ is hereditarily undecidable, making $T$ undecidable as required. \bs

One can see the key step in this corollary was using the following property inherent to the considered fields $L$:
\begin{align*}
\mbox{$K$ relatively algebraically closed in $L$} \implies K \equiv L. \tag{$\flat$}
\end{align*}

In $\S \ref{scf}$ we will outline Ziegler's construction of the field $K_q$, with a minor discrepancy for $L$ a separably closed field. In $\S \ref{hvfsect}$ we outline Ziegler's construction in the case $L=\Q_p$ but again with minor changes so his method works for a general class of henselian valued fields. Later in the paper we will discuss extending this construction to more difficult cases that avoid property $(\flat)$.

\bigskip
\section{Separably Closed Fields}\label{scf}

To save referring the reader to another text, we will outline Ziegler's construction in this subsection. To make this a more interesting exercise we shall {prove \emph{the theory of any separably closed field is finitely undecidable}; not considered by Ziegler in \cite{ziegler}. }

{Let {\ttfamily SCF} denote the theory of separably closed fields, {\ttfamily SCF}${}_{p}$ the theory of separably closed fields of characteristic $p$, and {\ttfamily SCF}${}_{p, e}$ the theory of separably closed fields of characteristic $p$ and degree of imperfection $e$. We shall assume for this section the reader is familiar with \cite{bouscaren}. As {\ttfamily SCF}${}_{0}$ $=$ {\ttfamily ACF}${}_0$ and {\ttfamily SCF}${}_{p, 0}$ $=$ {\ttfamily ACF}${}_p$, we will not consider {the cases $p = 0$ or $p > 0$ \& $e = 0$} (the case $e = \infty$ will be considered separately in \emph{Corollary \ref{scf2}}).}

Let $q \neq p$ be a prime number, {$L = \left(\widetilde{\mathbb{F}_p(t)}(u_1, \dots, u_e)\right)^s$ where $\{u_1, \dots, u_e\}$ are transcendental and algebraically independent over $\widetilde{\mathbb{F}_p(t)}$ \phantom{$L = \left(\widetilde{\mathbb{F}_p(t)}(u_1, \dots, u_e)\right)^s$} \hspace{-45mm} (and `$\,F^s\,$' denotes the \textit{separable} closure of a field $F$). Note $L \models \mbox{\ttfamily SCF}_{p,e}$, as $\widetilde{\mathbb{F}_p(t)}$ is perfect. First we have:}

\begin{prop}\label{sect1mainthm}
For each prime $q \neq p$, there exists a field $K_q \subseteq L$ such that:
\begin{enumerate}
    \item {$\mathbb{F}_p[t]$ is definable (with parameters)} in $K_q$;
    \item If the intermediate field $K_q \subseteq H \subseteq L$ is finite \emph{and separable} over $K_q$, then $[H : K_q] = 1$ or is divisible by $q$.
\end{enumerate}
\end{prop}

{To prove this, we require a construction. Let $F = \left(\widetilde{\mathbb{F}_p(t)}(u_1, \dots, u_{e-1})\right)^s \subseteq L$; we construct a field $K_q \subseteq L$ such that}
\begin{align*}
    F &= \{ a \in K_q \mbox{ : }\forall b \in K_q^* \mbox{ } ([1 + b \in (K_q)^q \land a^q + b^{-1} \in (K_q)^q] \rightarrow b \in (K_q)^q)\},\\
    \mathbb{F}_p[t] &= \{r \in F \mbox{ : } \forall r_1 \neq r_2 \in F \mbox{ s.t. } r_1 + r_2 = r\mbox{, } u_e^q - r_1 \mbox{ or } u_e^q - r_2 \in (K_q)^q\}.
\end{align*}

This will suffice to prove \emph{Proposition \ref{sect1mainthm} (1)}. We will take $K_q$ to be the union of a sequence
$$F(u_e) = E_0 \subseteq E_1 \subseteq E_2 \subseteq \dots$$

within $L$ of finite separable extensions $E_i/F(u_e)$. Obtaining \emph{Proposition \ref{sect1mainthm} (2)} while ensuring $F$ and $\mathbb{F}_p[t]$ are definable in this way requires us to keep a tight rein on the $q$-th roots in $K_q$. To that end, we will also carefully construct a sequence:
$$\emptyset = S_0 \subseteq S_1 \subseteq S_2 \subseteq \dots$$

of finite subsets $S_i \subseteq E_i$, ultimately desiring $K_q\setminus (K_q)^q = \bigcup_i S_i$. To ensure we do not introduce an incompatibility between $K_q \setminus (K_q)^q$ and $F$ or $\mathbb{F}_p[t]$, we will ask the following rule {(\cite[p.~273]{ziegler})} to be obeyed at each point of the sequence $(E_i, S_i)$:\label{5}

{\begin{align*}
  &\mbox{There is a family of valuations $\{v_s\}_{s \in S_i}$ on $E_i$ such that $v_s(F) = 0$ and}\\
  &\mbox{$q \nmid v_s(s)$ for $s \in S_i$. In addition, for all $r_1 \neq r_2 \in F$ with $r_1 + r_2 \in \mathbb{F}_p[t]$,}\hspace{7mm} \tag{$\clubsuit$}\\
  &\mbox{either $\forall s \in S_i$, $q | v_s(u_e^q - r_1)$, or $\forall s \in S_i$, $q | v_s(u_e^q - r_2)$.}\\
\end{align*}}

\noindent We will reference the following two standard lemmas, taken directly from \cite[\S 3]{ziegler}:

\begin{lem}\label{zieg1}
\emph{\cite[Lemma 1]{ziegler}.} Let $(H_1, v_1)$ be a discretely valued field, $H_2/H_1$ a finite extension, and $q$ a prime such that $q \nmid [H_2 : H_1]$. Then there is an extension of $v_1$ to $H_2$, which we denote $v_2$, such that $q \nmid (v_2 H_2 : v_1 H_1)$. \bsn
\end{lem}

\begin{lem}\label{zieg2}
\emph{\cite[Lemma 2]{ziegler}.} Let $(H, v)$ be a valued field and $q$ a prime distinct to $\charac(Hv)$. For $a \in H \setminus (H)^q$ with $q | v(a)$, there is an extension of valued fields $(H(\sqrt[q]{a}), w)/(H, v)$ such that $wH(\sqrt[q]{a}) = vH$. \bsn
\end{lem}

\noindent We will also require the following fact:

\begin{definition}\label{dependent}
Valuations $v_1, v_2$ on a field $K$ are \emph{dependent} if $\mathcal{O}_{v_1} \mathcal{O}_{v_2} \subsetneq K$.
\end{definition}

\begin{lem}\label{same}
If $v_1, v_2$ are dependent discrete valuations on a field $K$, then $v_1 = v_2$ $($by which we mean $\mathcal{O}_{v_1} = \mathcal{O}_{v_2})$.
\end{lem}

\proof
The valuation ring of a discrete valuation is maximal \cite[Corollary 2.3.2]{EP}, hence $\mathcal{O}_{v_1} = \mathcal{O}_{v_1}\mathcal{O}_{v_2} = \mathcal{O}_{v_2}$ as desired. \bs

The construction begins with an enumeration {$a_0, a_1, a_2, \dots$ of the elements of $L$} \underline{separably algebraic} over $F(u_e)$, each repeated countably infinitely many times. Suppose $(E_i, S_i)$ is already constructed -- Ziegler considers now four cases, based on the equivalence class of $i$ mod $4$. These correspond to guaranteeing \emph{Proposition \ref{sect1mainthm} (1)} (\textsc{Case 1}), \emph{Proposition \ref{sect1mainthm} (2)} (\textsc{Case 2}), the definition of $F$ by ensuring there is a `reason' $b$ is excluded from $(K_q)^q$ (\textsc{Case 3} + $\clubsuit$), and the definition of $\mathbb{F}_p[t]$ (\textsc{Case 4} + $\clubsuit$); again, by ensuring there is a `reason' $a_n$ is excluded. 

\bigskip
\noindent\textsc{Construction.} (cf.\ \cite[\S 3]{ziegler})\label{constructionbegin}
\begin{case}
    \item $i = 4n$. If $q | [E_i(a_n) : E_i]$, then $(E_{i+1}, S_{i+1})= (E_i, S_i)$. Otherwise set $(E_{i+1}, S_{i+1})$ $= (E_i(a_n), S_i)$ and using \emph{Lemma \ref{zieg1}} extend each valuation $v_s$, $s \in S_i$, from $E_i$ to $E_{i+1}$ in a way preserving $(\clubsuit)$.\\
    
    \item $i = 4n + 1$.\label{case2} {Unless $a_n \in E_i \setminus S_i$, set $(E_{i+1}, S_{i+1})= (E_i, S_i)$.} Otherwise, if for some $v_s$, $s \in S_i$, we have $q \nmid v_s(a_n)$ then define $(E_{i+1}, S_{i+1})= (E_i, S_i \cup \{a_n\})$ and set $v_{a_n} \defeq v_s$. This ensures $(\clubsuit)$ holds for $i+1$. If $q| v_s(a_n)$ for all $s \in S_i$, then we take $(E_{i+1}, S_{i+1}) =\!\!\footnotemark\,\, (E_i(\sqrt[q]{a_n}), S_i)$ \footnotetext{Note $L = (L)^q$ as $L$ is separably closed and $q \neq p$. In addition, we consider $\sqrt[q]{a_n} \in E_i$ if $a_n \in (E_i)^q$.}and extend every valuation according to \emph{Lemma \ref{zieg2}}.\\
    
    \item $i = 4n +2$. Unless $a_n \in E_i \setminus F$, let $(E_{i+1}, S_{i+1})= (E_i, S_i)$. If $a_n \in E_i \setminus F$ let $v$ be a discrete valuation on $E_i$, trivial on $F$, which is negative on $a_n$.     
    If the second condition in $(\clubsuit)$ does not already hold for $\{v, \{v_s\}_{s \in S_i}\}$ in $E_i$, then there exists $r \in F$ such that $q \nmid v(u_e^q - r)$ and $q | v_s(u_e^q - r)$ for all $s \in S_i$. By the strong triangle inequality, there is at most one such $r$: indeed, for $r \neq r' \in F$, $v(u_e^q - r') = v(u_e^q - r + r - r') = 0$, as $v(r - r') = 0$ and $q \nmid v(u_e^q - r)$. As $L = (L)^q$, we may set $E = E_i(\sqrt[q]{u_e^q - r})$ and extend the valuations $\{v, \{v_s\}_{s \in S_i}\}$ sensibly as above. We conclude the second condition of $(\clubsuit)$ holds for $(E, \{v, \{v_s\}_{s \in S_i}\})$. 
    
    {If $v$ is independent to $v_s$ for every $s \in S_i$: let $\{v, v_{s_1}, \dots, v_{s_k}\}$ be the distinct valuations of $\{v, \{v_s\}_{s \in S_i}\}$. By the \emph{Approximation Theorem} {\cite[Theorem 2.4.1]{EP}}, there exists $b \in E$ such that $v(b)$ is the smallest positive element in the value group of $v$ (hence $q | v(1+b)$, $q|v(a_n^q + b^{-1})$) and $q|v_{s_j}(b)$, $q|v_{s_j}(1+b)$, $q|v_{s_j}(a_n^q + b^{-1})$ for $1 \leq j \leq k$. As $b, 1+b, a_n^q + b^{-1} \in (L)^q = L$, we may define
    $$(E_{i+1}, S_{i+1}) = \left(E\!\left(\sqrt[q]{1+b}, \sqrt[q]{a_n^q + b^{-1}}\right)\!,\, S_i \cup \{b\}\right).$$    
    Extending $\{v_b = v, \{v_s\}_{s \in S_i}\}$ as above, we know $(\clubsuit)$ holds as it did on $E$. 
    
    If $v$ is dependent with $v_{\widehat{s}}$ for some $\widehat{s} \in S_i$: by \emph{Lemma \ref{same}} $v = v_{\widehat{s}}$. Let $\{v_{s_1}, \dots, v_{s_l}\}$ be the distinct valuations of $\{v, \{v_s\}_{s \in S_i}\}$, assuming WLOG $v = v_{s_1}$. By the \emph{Approximation Theorem} {\cite[Theorem 2.4.1]{EP}}, there exists $b \in E$ such that $v_{s_1}(b)$ is the smallest positive element in the value group of $v_{s_1}$, and $q|v_{s_j}(b)$, $q|v_{s_j}(1+b)$, $q|v_{s_j}(a_n^q + b^{-1})$ for $2 \leq j \leq l$. As $b, 1+b, a_n^q + b^{-1} \in (L)^q = L$, we may define
    $$(E_{i+1}, S_{i+1}) = \left(E\!\left(\sqrt[q]{1+b}, \sqrt[q]{a_n^q + b^{-1}}\right)\!,\, S_i \cup \{b\}\right).$$    
    Extending $\{v_b = v_{s_1}, \{v_s\}_{s \in S_i}\}$ as above, again $(\clubsuit)$ holds on $(E_{i+1}, S_{i+1})$. (This case allows $b \in S_i$ without issue, by \emph{Lemma \ref{same}}.)}\\
    
    \item $i = 4n+3$.\label{case4} Unless $a_n \in F\setminus \mathbb{F}_p[t]$ we {set $(E_{i+1}, S_{i+1}) = (E_{i}, S_{i})$.} Otherwise, first observe $B = \{r \in F \mbox{ : } \exists s \in S_i \mbox{ s.t. } q \nmid v_s(u_e^q - r)\}$ is finite. Next, for $r \in F^*$ there exists a discrete valuation $v_r$ on $F(u_e)$, trivial on $F$, for which $v_r(u_e^q - r)$ is the smallest positive element of its valuation group (this follows as $X^q - r \in F[X]$ has no multiple factors). For each such $r$, we choose an extension $w_r$ of $v_r$ to $E_i$, and by the construction of $E_i/F$, the set $C =$\linebreak $\{r \in F^* \mbox{ : } q|w_r(u_e^q - r)\}$ is finite. Choose $r_1 \in F^*$ such that $r_1 \neq a_n$, {$2 r_1 \neq a_n$,} and $r_1, a_n -r_1 \not\in C$, and $r_1, a_n - r_1 \not\in \mathbb{F}_p[t]$. Let $r_2 = a_n - r_1$ and finally define:
    $$(E_{i+1}, S_{i+1}) = (E_i, S_i \cup \{ u_e^q - r_1, u_e^q - r_2\}).$$
    One can prove $\{w_{r_1}, w_{r_2}, \{v_s\}_{s \in S_i}\}$ satisfies $(\clubsuit)$ based on this construction.\label{constructionend}\\
\end{case}

\begin{lem}\label{theprop}
Set $K_q = \bigcup_i E_i$. The above construction ensures we have the following features of $K_q$, and the definitions of $F$ and $\mathbb{F}_p[t]$ we intended:
\begin{enumerate}
    \item $F \subseteq (K_q)^q$.
    \item $K_q \setminus (K_q)^q = \bigcup_i S_i$.
    \item $F = \{a \in K_q \mbox{ \emph{:} } \forall b \in K_q^* \mbox{ } [( 1 + b \in (K_q)^q \land a^q + b^{-1} \in (K_q)^q) \rightarrow b \in (K_q)^q]\}$.
    \item $\mathbb{F}_p[t]\! =\! \{r \in F\! \mbox{ \emph{:} }\! \forall r_1 \neq r_2 \in F\! \mbox{ } (r_1 + r_2 = r)\! \rightarrow\! (u_e^q - r_1\! \in\! (K_q)^q \lor u_e^q - r_2\! \in\! (K_q)^q)\}$.
\end{enumerate}
\end{lem}

\proof
We follow \cite[\S 4]{ziegler} as much as possible.
\begin{enumerate}
    \item As $F$ is separably closed, and $q \neq p$, $F = (F)^q$.\\
    
    \item {Let $a \in (K_q)^q$. For all sufficiently large $i$, $a \in (E_i)^q$; hence $q | v(a)$ for all $v$ trivial on $F$. Therefore by $(\clubsuit)$ we have $a \not\in S_i$; consequently $a \not\in \bigcup_i S_i$. Conversely if $a \in K_q\setminus (K_q)^q$, then for some $n$ sufficiently large we have $a = a_n$ and $a \in E_{4n+1}$. By {\textsc{Case 2}} of the construction, $a \in S_{4n+2}$. This proves $K_q \setminus (K_q)^q = \bigcup_i S_i$.}\\
    
    \item Fix $a \in F$. Suppose for some nonzero $b \in K_q$ that $1 + b, a^q + b^{-1} \in (K_q)^q$. Let $i$ be so large that $1+b, a^q + b^{-1} \in (E_i)^q$. Let $v$ be any valuation on $E_i$ that is trivial on $F$. If $v(b) > 0$, then $v(b) = -v(a^q + b^{-1})$ is divisible by $q$. If $v(b) < 0$, then $v(b) = v(1+b)$ is divisible by $q$. Hence $q| v(b)$ always. {By $(\clubsuit)$, $b \not\in S_i$; by (2) therefore $b \in (K_q)^q$}. Conversely, if $a \in K_q \setminus F$, we may choose $n$ sufficiently large such that $a = a_n \in E_{4n+2}$. {In {\textsc{Case 3}}} we make it such that in $S_{4n+3}$ there is a (nonzero) $b$ with $1 + b, a^q + b^{-1} \in (E_{4n+3})^q$. This concludes the proof.\\
    
    \item Let $r_1 + r_2 \in \mathbb{F}_p[t]$, with $r_1 \neq r_2 \in F$. If it is the case that $u_e^q - r_1, u_e^q - r_2 \not\in (K_q)^q$, then for some sufficiently large $i$, they belong to $S_i$. However this contradicts $(\clubsuit)$. If we suppose $r \in F \setminus \mathbb{F}_p[t]$, for some sufficiently large $n$ it is the case that $a_n = r$. Then by {\textsc{Case 4}} there exists $r_1 \neq r_2 \in F$, $r_1 + r_2 = r$, such that $u_e^q - r_1, u_e^q - r_2 \in S_{4n+4}$. By {(2)} this ensures $u_e^q - r_1, u_e^q - r_2 \not\in (K_q)^q$; {again} a contradiction. \bs
\end{enumerate}

\noindent \textit{Proof of \emph{Proposition \ref{sect1mainthm}:}}\\
First, as $F$ and $\mathbb{F}_p[t]$ are definable, $\N$ is interpretable (with parameters) as an $\mathcal{L}_r$-structure in $K_q$. Next, note that $K_q/F(u_e)$ is a separable extension, as by construction it is a union of finite separable extensions, and $K_q \subseteq F(u_e)^s$. Let $K_q \subset H \subseteq L$ be a finite separable extension. Then $H = K_q(a)$ for some $a \in L$ by the \textit{Primitive Element Theorem}. As $K_q(a)/F(u_e)$ is separable, for some $n$ sufficiently large we have $a = a_{n}$ and 
$$[E_{4n}(a_n) : E_{4n}] = [K_q(a) : K_q],$$
as we assume $a \not\in K_q$. By construction, $q | [E_{4n}(a) : E_{4n}]$. \bs

Combining these fields in a nonprincipal ultraproduct, as will be done in the next corollary, allows us to conclude the desired undecidability result.

\begin{cor}\label{scf1}
Let $p$ be a prime and $e \in \N_{>0}$. Then \emph{\ttfamily SCF}${}_{p, e}$ is finitely undecidable.
\end{cor}

\proof
Let $L = \left(\widetilde{\mathbb{F}_p(t)}(u_1, \dots, u_e)\right)^s$ as before \emph{Proposition \ref{sect1mainthm}}. Let $P$ be the set of primes distinct to $p$. For each $q \in P$, use \emph{Proposition \ref{sect1mainthm}} to obtain a field $K_q$ satisfying \emph{(1) \& (2) ibid.} Let $\mathcal{U}$ be a nonprincipal ultrafilter on $P$ and let $\mathbb{K}$ be the ultraproduct $\prod_{q \in P} K_q / \mathcal{U}$. We claim that $\mathbb{K}$ is relatively separably closed in $L^{\mathcal{U}}$: indeed, suppose $f = (f_q) \in \mathbb{K}[X]$ is separable and has a root $\alpha = (\alpha_q) \in L^{\mathcal{U}}$. {In which case
\begin{align*}
  &\{q \mbox{ : } f_q(x) \in K_q[X] \mbox{ separable, and } f_q(\alpha_q) = 0\} \in \mathcal{U},\\
\mbox{hence }&\{q \mbox{ : } K_q(\alpha_q)/K_q \mbox{ separable and } [K_q(\alpha_q) : K_q] \leq \deg(f)\} \in \mathcal{U}.  
\end{align*}
By \emph{Proposition \ref{sect1mainthm} (2)}, $\{q \mbox{ : } [K_q(\alpha_a) : K_q] = 1\} \in \mathcal{U}$, and thus $\alpha \in \mathbb{K}$ as desired.}

Therefore $\mathbb{K}$ is a separably closed field of characteristic $p$. Recall $\{u_1, \dots, u_e\}$ is a $p$-basis for $L$. As they are $p$-independent in $L$, and by construction $u_1, \dots, u_{e} \in K_q$, they remain $p$-independent in $K_q$. Hence the degree of imperfection of $K_q$ is at least $e$, for each $q$. Moreover, by construction $K_q /\widetilde{\mathbb{F}_p(t)}(u_1, \dots, u_e)$ is an algebraic (separable) extension. As algebraic extensions do not increase the degree of imperfection, the degree of imperfection of $K_q$ is at most $e$. Thus by \textit{\L o\'{s}' Theorem}, the degree of imperfection of $\mathbb{K}$ is exactly $e$. We conclude that $\mathbb{K} \models \mbox{\ttfamily SCF}_{p, e}$, and hence $\mathbb{K} \models T$ for any finite subtheory $T \subseteq \mbox{\ttfamily SCF}_{p, e}$.

Since $T$ is finitely axiomatised, there exists some prime $q$ such that $K_q \models T$. {By \emph{Lemma \ref{theprop} \& Corollary \ref{appendixacor}}, $\Th(K_q; \mathcal{L}_r)$ is hereditarily undecidable, making $T$ undecidable as required.} \bs

\begin{cor}\label{scf2}
Let $p$ be prime. Then \emph{\ttfamily SCF}${}_{p, \infty}$ is finitely undecidable.
\end{cor}

\proof
{Let $T$ be a finite subtheory of {\ttfamily SCF}${}_{p, \infty}$, which we assume is axiomatised by the axioms of a field of characteristic $p$, such that each separable polynomial over the field has a root in the field, and for each $n \in \N_{>0}$ the statement ``the degree of imperfection is greater than $n$''. By the \emph{Compactness Theorem}, there exists a finite subset $\Delta$ of this axiomatisation such that $\Delta \models T$. For some finite $\nu$ sufficiently large, {\ttfamily SCF}${}_{p, \nu} \models \Delta$, hence $T$ is a finite subtheory of {\ttfamily SCF}${}_{p, \nu}$. The result follows from \emph{Corollary \ref{scf1}}.} \bs

\begin{exmp}
For all primes $p > 0$ and $e \in \N \cup \,\{\infty\}$, the theory $\mbox{\ttfamily SCF}_{p,e}$ is known to be decidable (see \cite[pp.\! 146--153]{bouscaren} for exposition). Therefore \emph{Corollaries \ref{scf1} \&\ \ref{scf2}} put a bound on further possible decidability results for these theories. \sq
\end{exmp}

It is worth remarking that, modulo some conjectures, these results are in connection with aspects of classification theory. It is a theorem of Macintyre \cite{macintyre} that every infinite $\omega$-stable field is a model of {\ttfamily ACF}${}_p$ for $p = 0$ or prime. From the 1970's we have the following conjecture:

\begin{customconj}{(Stable Fields)}
\textit{Every infinite stable field is separably closed.} \bsn
\end{customconj}

This is known in some cases, such as for the aforementioned $\omega$-stable \cite{macintyre} or superstable infinite fields \cite{cherlin-shelah}, or for infinite stable fields of weight 1 \cite{krupinsky-pillay} or finite dp-rank \cite{halevi-palacin}, or most recently infinite large stable fields \cite{jtwy}.

\begin{cor}
Assume the Stable Fields Conjecture. Then every infinite stable field is finitely undecidable. \bsn
\end{cor}

Let us use this connection to classification theory to motivate which fields to consider next. Outside of stable theories, there are two orthogonal directions in which to travel: one direction attempts to extend the theories of forking, dividing and independence of types to more general contexts (e.g.\ [super]simple, [super]rosy), while the other direction aims to understand theories with a modest notion of order (e.g.\ o-minimal, NIP). The latter direction contains theories we are already familiar with: {\ttfamily RCF} in the language of ordered rings is o-minimal, \textit{p}{\ttfamily CF} in the language of valued fields is distal and dp-minimal (hence NIP). One might wonder what other field theories could be present under this banner -- and there is a conjecture of Shelah that would answer this question:

\begin{customconj}{(Shelah/NIP Fields)}
\textit{Every infinite NIP field is either separably closed, real closed, or admits a nontrivial henselian valuation.}
\end{customconj}

\begin{thm}\label{maybejochen}
Assume the NIP Fields Conjecture. Then every infinite NIP field is either real closed, separably closed, or admits a nontrivial henselian valuation \emph{$\emptyset$-definable in the language of rings}. 
\end{thm}

\proof
(Here a valuation is \textit{definable} if the valuation ring is a definable subset of the field.) $\mathcal{L}_r$-definability is \cite[Proposition 6.2(2)]{hhj}, and the results cited in the proof (from \cite{jochenjahnke}) in fact conclude $\emptyset$-definability. \bs

Therefore a sensible goal would be to prove that every field with a nontrivial $\emptyset$-$\mathcal{L}_r$-definable henselian valuation is finitely undecidable. Or more so, that every henselian valued field is finitely undecidable in the language of valued fields $\mathcal{L}_{val}$.

\bigskip
\section{Equicharacteristic 0 Henselian Valued Fields}\label{hvfsect}

The previous subsection did not address the aspects of Ziegler's construction relevant to $\Q_p$; these aspects will be seen in this subsection. In this subsection we will consider a pair of valued fields $(R, v_R)$, $(Z, v_Z)$, and an additional field $F$, under the following assumptions:

\begin{customass}{($\otimes$)}
\phantom{boo}
\emph{\begin{enumerate}
    \item $R \subseteq F \subseteq Z$, $v_R = v_Z|_R$, and $(Z, v_Z)$ is a henselian immediate extension of $(R, v_R)$;
    \item $R$ (thence $v_R R$ and $R v_R$) is countable, and if $\charac(R) > 0$ then $R$ is transcendental over its prime subfield;
    \item There are uncountably many elements of $Z$ transcendental over $R$;
    \item $F = Z \cap {R(\overline{x})}^s$, where $R(\overline{x})$ is a purely transcendental, finite transcendence \linebreak degree extension of $R$;
    \item Let $q > \charac(Rv_R)$ be prime; then $Z = (Z)^q \cdot F^*$.
\end{enumerate}}
\end{customass}

First we will give a concrete example of a pair of valued fields where these assumptions are satisfied. Let $k$ be a field and $\Gamma$ an ordered abelian group. Consider the multiplicative group of formal monomials $\{t^{\gamma} \mbox{ : } \gamma \in \Gamma\}$, where $t^0 = 1$ and $t^{\gamma_1} \cdot t^{\gamma_2} = t^{\gamma_1 + \gamma_2}$. Define $k[\Gamma]$ to be the set of formal series $\sum_{\gamma} a_{\gamma} t^{\gamma}$ where $a_{\gamma} \in k$ and only finitely many $a_{\gamma}$ are nonzero. Addition and multiplication are defined by

\begin{align*}
        \sum_{\gamma} a_{\gamma} t^{\gamma} + \sum_{\gamma} b_{\gamma} t^{\gamma} &= \sum_{\gamma} (a_{\gamma} + b_{\gamma}) t^{\gamma};\\
        \left( \sum_{\gamma} a_{\gamma} t^{\gamma} \right) \cdot \left( \sum_{\gamma} b_{\gamma} t^{\gamma} \right) &= \sum_{\gamma} \left(\sum_{\gamma_1 + \gamma_2 = \gamma} a_{\gamma_1} b_{\gamma_2} \right) t^{\gamma}.\\
\end{align*}
    
These operations are confirmed to be well-defined, and $k[\Gamma]$ an integral domain, by \cite[\S 2.4]{markervf}. This domain comes with a natural valuation $v_{\Gamma}(\sum_{\gamma} a_{\gamma} t^{\gamma}) \defeq \min \supp(\sum_{\gamma} a_{\gamma} t^{\gamma})$. Define $k(\Gamma)$ to be the fraction field of this valued domain. Further define $k((\Gamma))$ as the set whose elements are formal series $\sum_{\gamma} a_{\gamma} t^{\gamma}$ with well-ordered support. By \cite[\S 2.4]{markervf}, $(k((\Gamma)), v_{\Gamma})$ is a well-defined immediate henselian overfield of $(k(\Gamma), v_{\Gamma})$.

\begin{lem}\label{concrete}
Let $e \in \N$, $k$ and $\Gamma$ be countable, and $v$ be a henselian valuation on $k((\Gamma))$ which factors through $v_{\Gamma}$, i.e.\ there exists a valuation $v'$ on $k$ such that\footnote{By ``$\,v' \circ v_{\Gamma}$'' we mean the composition of places $res_{v'}$, $res_{v_{\Gamma}}$, which give rise to the valuations as per usual (cf.\ \cite[Construction 2.2.6]{friedjarden}).} $v = v' \circ v_{\Gamma}$. There exists $t_1, \dots, t_e \in k((\Gamma))$ transcendental over $k(\Gamma)$ and algebraically independent such that the pair of valued fields $(R, v_R) = (k(\Gamma), v|_{k(\Gamma)})$, $(Z, v_Z) = (k((\Gamma)), v)$, and $F = k((\Gamma)) \cap \left(k(\Gamma)(t_1, \dots, t_e)\right)^s$ satisfy \emph{Assumption $(\otimes)$}.
\end{lem}

\proof
Properties \emph{(1) \& (4)} follow by definition. Property \emph{(2)} follows by construction, and as $k(\Gamma)$ is countable (from its definition). Property \emph{(3)} can be seen by a cardinality argument (cf.\ \cite[p.~82]{vdd}): fixing $\gamma \in \Gamma^{>0}$, there is an injection $(\N; 0, +, <) \hookrightarrow (\Gamma; 0, +, <)$ given by $n \mapsto n \cdot \gamma$, and by definition $|k((\Gamma))| \geq |k|^{\aleph_0} \cdot 2^{\aleph_0} = 2^{\aleph_0}$, while $|\widetilde{k(\Gamma)}| = \aleph_0$.

Property \emph{(5)} requires more work: we adapt \cite[Lemma 3]{ziegler}. Clearly $(k((\Gamma)))^q \cdot F^* \subseteq k((\Gamma))$; we are required to show that for all $a \in k((\Gamma))^*$, there exists $b \in F^*$ such that $a b^{-1} \in (k((\Gamma)))^q$. Choose $b \in F^*$ such that $v_{\Gamma}(a - b) > v_{\Gamma}(a)$; this can be done by setting $b \in k(\Gamma)^* \subseteq F^*$ to be a sufficiently large finite truncation of $a$. Then $v_{\Gamma}(a b^{-1} - 1) > v_{\Gamma}(a b^{-1})$, hence $ab^{-1} \equiv 1$ mod $\mathfrak{m}_{v_{\Gamma}}$. By Hensel's Lemma (regardless of $v$, $(k((\Gamma)), v_{\Gamma})$ is henselian), $ab^{-1}$ is a $q$-th power in $k((\Gamma))$, as desired. \bs

Using \emph{Assumption $(\otimes)$}, for $q > p$ prime, we shall construct a field extension $F \subseteq K_q \subseteq Z$ such that $\Z$ or $\mathbb{F}_p[z]$ (where $p = \charac(R) > 0$ and $z \in Z$ is transcendental over $\mathbb{F}_p$) is $\mathcal{L}_{val}$-definable in $K_q$, and for \underline{some elements} $a \in Z$ algebraic over $K_q$, $q|[K_q(a) : K_q]$. 

\begin{remark}\label{littleprob}
Notice that if $Z$ is perfect, and for all $a \in Z$ algebraic over $K_q$ we have either $K_q(a) = K_q$ or $q | [K_q(a) : K_q]$, then $K_q$ is perfect: $(K_q)^p = K_q$. This will be a problem for \emph{Theorem \ref{nipequichar}}, where we will consider finite subtheories $T$ of \textit{imperfect} fields and prove $K_q \models T$. This problem will be resolved after \emph{Lemma \ref{independenceyay}}. \sq
\end{remark}

By \emph{Assumption $(\otimes)$} there exists an element $t \in Z$ transcendental over $F$. The field $K_q$ will be the union of a specific sequence of finite extensions of $F(t)$ in $Z$:
$$R \subseteq F \subset F(t) =  E_0 \subseteq E_1 \subseteq E_2 \subseteq \cdots \quad \subseteq Z.$$

As before, we will also construct a sequence $\emptyset = S_0 \subseteq S_1 \subseteq S_2 \subseteq \cdots$ of finite subsets $S_i \subseteq E_i \cap (Z)^q$, ultimately desiring a close relationship between $K_q \setminus (F^* \cdot (K_q)^q)$, $(K_q \cap (Z)^q) \setminus (K_q)^q$, and $\bigcup_i S_i$. We will desire $F$ \&\ $\Z$ (if $\charac(R) = 0$; $\mathbb{F}_p[z]$ otherwise) to have the following definitions, similar to as before:

\begin{align*}
    &F = \{ a \in K_q \mbox{ : }\forall b \in K_q \cap (Z)^q\setminus\{0\} \mbox{, } ([1 + b \in (K_q)^q \land a^q + b^{-1} \in (K_q)^q] \rightarrow b \in (K_q)^q)\},\\
    &\Z \mbox{ (resp.\! $\mathbb{F}_p[z]$)} = \{u \in F \mbox{ : } \forall u_1 \neq u_2 \in F \mbox{ s.t. } u_1 + u_2 = u\mbox{, } t^q - u_1 \mbox{ or } t^q - u_2 \in (K_q)^q\}.\\
\end{align*}

Denote by ``$v_Z|_{E_i}$'' the restriction of $v_Z$ to $E_i \subseteq Z$. To again ensure we do not introduce an incompatibility between $(K_q \cap (Z)^q) \setminus (K_q)^q$ and $F$ or $\Z$ (resp.\! $\mathbb{F}_p[z]$), the following rule will be enforced during the construction:
\pagebreak 

\begin{align*}
  &\mbox{There is a family of discrete valuations $\{v_s\}_{s \in S_i}$ on $E_i$ such that $v_s(F) = 0$}\\
  &\mbox{and $q \nmid v_s(s)$ for $s \in S_i$. In addition, for all $u_1 \neq u_2 \in F$ with $u_1 + u_2 \in \Z$}\hspace{7mm} \tag{$\raisebox{-0.5mm}{$\varheart$}$}\\
  &\mbox{(resp.\! $\mathbb{F}_p[z]$), either $\forall s \in S_i$, $q | v_s(t^q - u_1)$, or $\forall s \in S_i$, $q | v_s(t^q - u_2)$.}\\
\end{align*}
We have the following lemma:

\begin{lem}\label{independenceyay}
Let $u$ be a nontrivial discrete valuation on $E_i$, considered as a (finite) function field extension of $F(t)/F$. Then $u$ and $v_Z|_{E_i}$ are independent $($in the sense of \emph{Definition \ref{dependent})}.
\end{lem}

\proof
Assume $u$ and $v_Z|_{E_i}$ are dependent: by \cite[Theorem 2.3.4]{EP} they induce the same topology on $E_i$. Thus there exists $a \in E_i$ such that $a \cdot \mathfrak{m}_{u} \subseteq \mathfrak{m}_{v_Z|_{E_i}}$ ($= \mathfrak{m}_{v_Z} \cap E_i$). However, as $R$ is a subset of the constant subfield of $E_i$ (and thus $u(R) = 0$, while $v_R R = v_Z Z$), there exists $f \in a \cdot \mathfrak{m}_u$ with $v_Z(f) = v_Z|_{E_i}(f) < 0$; a contradiction. \bs

To use the same construction of $K_q$ throughout the paper, yet have $K_q$ a model of a finite subtheory of {either} a \textit{perfect or imperfect} field, we have the following definition:

\begin{definition}
The set of \emph{pliant} elements over a field will denote either the \emph{separably algebraic} \textbf{or} \emph{algebraic} elements over the field, and be specified in \emph{Corollary \ref{thisis}/Theorem \ref{newmixed}/Theorem \ref{nipequichar}.}
\end{definition}

Fix an enumeration $a_0, a_1, \dots$ of the elements of $Z$ pliant over $F(t)$, each repeated countably infinitely many times. Suppose $(E_i, S_i)$ is already constructed -- and consider the following modified construction.

\bigskip
\noindent\textsc{Modified Construction.}\label{modifiedbegin}
\begin{case}
    \item $i = 4n$. As on p.\ \pageref{constructionbegin}.\\
    
    \item $i = 4n + 1$. If $a_n \not\in E_i$ or $a_n \not\in (Z)^q$ then set $(E_{i+1}, S_{i+1})= (E_i, S_i)$. Otherwise proceed as on p.\ \pageref{case2}.\\ 
    
    \item $i = 4n+2$. Unless $a_n \in E_i \setminus F$, let $(E_{i+1}, S_{i+1})= (E_i, S_i)$. If $a_n \in E_i \setminus F$ let $w$ be a discrete valuation on $E_i$ (considered as a function field extension of $F(t)/F$), trivial on $F$, which is negative on $a_n$. By \emph{Lemma \ref{independenceyay}}, $w$ and $v_Z|_{E_i}$ are independent. Let us define a finite separable extension $E/E_i$: if the second condition of $($\raisebox{-0.5mm}{$\varheart$}$)$ already holds for $\{w, \{v_s\}_{s \in S_i}\}$ in $E_i$, set $E = E_i$. Otherwise there exists $u \in F$ such that $q \nmid w(t^q - u)$ and $q | v_s(t^q - u)$ for all $s \in S_i$. By the strong triangle inequality, there is at most one such $u$. 
    
    Under \emph{Assumption $(\otimes)(5)$} there exists $d \in F^*$ such that $d(t^q - u) \in (Z)^q$. Thus we may set $E = E_i(\sqrt[q]{d(t^q - u)})$ and extend the valuations $\{v_s\}_{s \in S_i}$ sensibly as before (and let $w'$ be any extension to $E$ of $w$). We conclude the second condition of $($\raisebox{-0.5mm}{$\varheart$}$)$ now holds for $(E, \{w', \{v_s\}_{s \in S_i}\})$.
    
    {If $w'$ is independent to $v_s$ for every $s \in S_i$: let $\{w', v_{s_1}, \dots, v_{s_k}\}$ be the distinct valuations of $\{w', \{v_s\}_{s \in S_i}\}$. By the \emph{Approximation Theorem} \cite[Theorem 2.4.1]{EP}, there exists $b \in E$ such that $w'(b)$ is the smallest positive element of its value group, $q | v_{s_j}(b)$ and $v_{s_j}(b) < 0, -v_{s_j}(a_n^q)$ for $1 \leq j \leq k$. Notice $q | w'(1+b)$, $q|w'(a_n^q + b^{-1})$, and $q|v_s(a_n^q + b^{-1})$, $q|v_s(1+b)$ for all $s \in S_i$. 
    
    We also wish $b, 1+b, a_n^q  + b^{-1} \in (Z)^q$. This can be achieved with further care by using the \emph{Approximation Theorem} and Hensel's Lemma: by \emph{Lemma \ref{independenceyay}} and \cite[Corollary 2.3.2]{EP}, $v_Z|_{E}$, $w'$, and $v_{s_j}$ for $1 \leq j \leq k$ are pairwise independent. Let $d \in (E)^q \subseteq (Z)^q$ have $v_Z|_E(d) > 0$. Using the \emph{Approximation Theorem}, we choose $b \in E$ so that in addition $v_Z|_E(b - d) > v_Z|_E(d)$. Then by Hensel's Lemma, $b \in E \cap (Z)^q$, and since $v_Z|_E((b+1) - 1) = v_Z|_E(b)$ and $v_Z|_E((b + a_n^{-q}) - a_n^{-q}) = v_Z|_E(b)$, we have $1 + b, b + a_n^{-q} \in E \cap (Z)^q$. (Hence $a_n^q + b^{-1} \in E \cap (Z)^q$.) We define
    $$(E_{i+1}, S_{i+1}) = (E(\sqrt[q]{\mbox{\normalsize$1+b$}}, \sqrt[q]{a_n^q  + b^{-1}}),\, S_i \cup \{b\}).$$    
    Extending $\{v_b = w', \{v_s\}_{s \in S_i}\}$ as in \textsc{Case 2}, $($\raisebox{-0.5mm}{$\varheart$}$)$ holds as it did on $E$.
    
    Now assume $w'$ is dependent with $v_{\widehat{s}}$ for some $\widehat{s} \in S_i$: in which case, $w' = v_{\widehat{s}}$ by \emph{Lemma \ref{same}}. Let $\{v_{s_1}, \dots, v_{s_l}\}$ be the distinct valuations of $\{w', \{v_s\}_{s \in S_i}\}$, assuming WLOG $w' = v_{s_1}$. By \emph{Lemma \ref{independenceyay}} and \cite[Corollary 2.3.2]{EP}, $v_Z|_{E}$, and $v_{s_j}$ for $1 \leq j \leq l$ are pairwise independent. Let $d \in (E)^q \subseteq (Z)^q$ have $v_Z|_E(d) > 0$. Using the \emph{Approximation Theorem}, we choose $b \in E$ so that $v_Z|_E(b - d) > v_Z|_E(d)$, $v_{s_1}(b)$ is the smallest positive element of its value group, $q | v_{s_j}(b)$ and $v_{s_j}(b) < 0, -v_{s_j}(a_n^q)$ for $2 \leq j \leq l$. By Hensel's Lemma, $1 + b, a_n^q + b^{-1} \in E \cap (Z)^q$, and we may define
    $$(E_{i+1}, S_{i+1}) = (E(\sqrt[q]{\mbox{\normalsize$1+b$}}, \sqrt[q]{a_n^q  + b^{-1}}),\, S_i \cup \{b\}).$$    
    Extending $\{v_b = w', \{v_s\}_{s \in S_i}\}$ as in \textsc{Case 2}, $($\raisebox{-0.5mm}{$\varheart$}$)$ holds on $(E_{i+1}, S_{i+1})$. (Again this case allows for $b \in S_i$ without issue, by \emph{Lemma \ref{same}}.)}\\
    
    \item $i = 4n +3$. As on p.\ \pageref{constructionend}. Recall this step extends $\{v_s\}_{s \in S_i}$ to $\{v_s\}_{s\in S_{i+1}} = \{w_{r_1}, w_{r_2}, \{v_s\}_{s \in S_i}\}$. The valuations $w_{r_1}, w_{r_2}$ are independent to $v_Z|_{E_i}$ by \emph{Lemma \ref{independenceyay}}, and $($\raisebox{-0.5mm}{$\varheart$}$)$ is satisfied. \label{modifiedend}\\
\end{case}

\noindent Now we can show the following (cf.\ \emph{Lemma \ref{theprop}}):

\begin{lem}\label{thelittlelem}
Set $K_q = \bigcup_i E_i$. The above construction ensures we have the following features of $K_q$ under \emph{Assumption $(\otimes)$:}
\begin{enumerate}
    \item $(K_q \cap (Z)^q) \setminus (K_q)^q = \bigcup_i S_i$.
    \item $F^* \cdot \left( \bigcup_i S_i \right) = K_q \setminus (F^* \cdot (K_q)^q)$.
    \item $F \!=\! \{ a \in {K_q}\!\mbox{ \emph{:} }\!\forall b \in K_q \cap (Z)^q\setminus\{0\} \mbox{ } ([1 + b \in (K_q)^q \land a^q + b^{-1} \in (K_q)^q]\! \rightarrow\! b \in (K_q)^q)\}$.
    \item $\Z \mbox{ $($resp.\! $\mathbb{F}_p[z])$} = \{u \in F \mbox{ \emph{:} } \forall u_1 \neq u_2 \in F \mbox{ s.t. } u_1 + u_2 = u\mbox{, } t^q - u_1 \mbox{ or } t^q - u_2 \in F^* \cdot (K_q)^q\}$.
\end{enumerate}
\end{lem}

\proof 
In \cite[\S 4]{ziegler}, but for exposition:
\begin{itemize}
    \item[(1)] Let $a \in (K_q \cap (Z)^q)\setminus (K_q)^q$. For some $n$ sufficiently large, $a_n = a$ and $a_n \in E_{4n+1}$. \textsc{Case 2} of the above construction assures $a_n \in S_{4n+2}$, hence $a \in \bigcup_i S_i$ as desired. Conversely, by construction $\bigcup_i S_i \subseteq K_q \cap (Z)^q$, and $S_i \cap (K_q)^q = \emptyset$ for all $i$.\\
    
    \item[(2)] Let $a \in F^* \cdot (K_q)^q$. For all $i$ sufficiently large, $a \in F^* \cdot (E_i)^q$ and hence $q|u(a)$ for all valuations $u$ trivial on $F$. By design of $($\raisebox{-0.5mm}{$\varheart$}$)$, $a \not\in F^* \cdot S_i$, hence $F^* \cdot \left( \bigcup_i S_i \right) \subseteq K_q \setminus (F^* \cdot (K_q)^q)$. Conversely, if $a \in K_q \setminus (F^* \cdot (K_q)^q)$, then by \emph{Assumption $(\otimes)(5)$} there exists $b \in F^*$ with $ab \in (K_q \cap (Z)^q) \setminus (K_q)^q = \bigcup_i S_i$ by (1).\\
    
    \item[(3)]  Let $a \in F$. Suppose for $b \in {K_q} \cap (Z)^q \setminus\{0\}$, we have $1 + b, a^q + b^{-1} \in (K_q)^q$. Let $i$ be sufficiently large such that $1 + b, a^q  + b^{-1} \in (E_i)^q$. Notice that, for any valuation $u$ on $E_i$ trivial on $F$, $q|u(b)$: indeed, if $u(b) < 0$ then $u(b) = u(1 + b)$, and if $u(b) > 0$ then $u(b) = u(a^q  + b^{-1})$. By $($\raisebox{-0.5mm}{$\varheart$}$)$ $b \not\in S_i$ (for all subsequent $i$ too), hence as $b \in K_q \cap (Z)^q$, $b \in (K_q)^q$ by (1). 
    
    Conversely, if $a \in K_q \setminus F$, then for some $n$ sufficiently large we may assume $a_n = a$ and $a \in E_{4n+2}$. By \textsc{Case 3} of the construction, there exists $b \in S_{4n+3}$ such that $1+ b, a^q + b^{-1} \in (E_{4n+3})^q$ and $b \in K_q \cap (Z)^q$. Therefore by (1),
    $$\exists b \in K_q \cap (Z)^q\setminus\{0\} \mbox{, }(1 + b \in (K_q)^q \land a^q + b^{-1} \in (K_q)^q \land b \not\in (K_q)^q),$$
    as desired.\\
    
    \item[(4)] Let $u \in \Z$ (resp.\! $\mathbb{F}_p[z]$), $u_1 \neq u_2 \in F$, and $u_1 + u_2  = u$. Assume for the purpose of contradiction both $t^q - u_1, t^q - u_2 \not\in F^* \cdot (K_q)^q$. By the argument in \textsc{Case 4} there exist $d_1, d_2$ such that $d_1(t^q - u_1), d_2(t^q - u_2) \in (K_q \cap (Z)^q) \setminus (K_q)^q = \bigcup_i S_i$, by (1). We conclude for some sufficiently large $i$ that $t^q - u_1, t^q - u_2 \in F^* \cdot S_i$, contradicting $($\raisebox{-0.5mm}{$\varheart$}$)$.
    
    Conversely, if $u \in F \setminus \Z$ (resp.\! $F \setminus \mathbb{F}_{p}[z]$), then for some $n$ sufficiently large we may assume $a_n = u$ and $a_n \in E_{4n+3}$. By \textsc{Case 3} of the construction, deliberately there exists $u_1 \neq u_2 \in F^*$ with $u_1 + u_2 = u$ and $t^q - u_1, t^q - u_2 \in S_{4n+4} \subseteq F^* \cdot S_{4n+4}$. Therefore by (2), $t^q - u_1, t^q - u_2 \not\in F^* \cdot (K_q)^q$ as required for this argument. \bs
\end{itemize}
\label{p14}

Let us return to the case $(R, v_R) = (k(\Gamma), v|_{k(\Gamma)})$, $(Z, v_Z) = (k((\Gamma)), v)$, $F = k((\Gamma)) \cap \left({k(\Gamma)(t_1, \dots, t_e)}\right)^s$. We have the following additional results:

\begin{thm}\label{akeref}
\emph{\textbf{(Ax-Kochen-Ershov)} \cite{axkochen, ershovake}.} Let $(K, v)$, $(L, w)$ be equicharacteristic $0$ henselian valued fields. Then $K \equiv_{\mathcal{L}_{val}} L$ if and only if $Kv \equiv_{\mathcal{L}_r} Lw$ and $vK \equiv_{\mathcal{L}_{oag}} wL$. 

Consequently, $(K, v) \equiv (Kv((vK)), v_{vK})$, the field of Hahn series in $vK$ over $Kv$.\bsn
\end{thm}

\begin{lem}\label{littledeflem}
Let $q > \charac(k)$ be prime and $v$ a henselian valuation on $k((\Gamma))$ which factors through $v_{\Gamma}$, where $k$ and $\Gamma$ are countable, and $K_q$ as above. Then $K_q \cap (k((\Gamma)))^q$ is $\mathcal{L}_{val}$-definable in $(K_q, w)$, where $w = v|_{K_q}$.
\end{lem}

\proof
Recall from \emph{Lemma \ref{concrete}} that \emph{Assumption $(\otimes)$} is satisfied. It suffices to show $c \in K_q \cap (k((\Gamma)))^q$ if and only if there exists $d \in K_q$ such that $w(c - d^q) > w(c)$.

Assume there exists $d \in K_q$ such that $w(c - d^q) > w(c)$; then (as elements of $k((\Gamma))$) we have $v(c - d^q) > v(c)$, hence $v(1 - \tfrac{d^q}{c}) > 0$, and thus $1 \equiv {\tfrac{d^q}{c}}$ mod $\mathfrak{m}_{v}$. By Hensel's Lemma, there exists $e \in k((\Gamma))$ such that $e^q = \tfrac{d^q}{c}$; we conclude $c \in (k((\Gamma)))^q$.

Conversely, let $c \in K_q \cap (k((\Gamma)))^q$ and write $c = \widehat{d}^q$. Let $d \in k(\Gamma)$ be a sufficiently large finite truncation of $\widehat{d} \in k((\Gamma))$ such that $v_{\Gamma}(\widehat{d} - d) > v_{\Gamma}(\widehat{d})$ (and note $v_{\Gamma}(\widehat{d}) = v_{\Gamma}(d)$). Then
\begin{align*}
v_{\Gamma}(\widehat{d}^q - d^q) &= v_{\Gamma}(\widehat{d} - d) + v_{\Gamma}(\widehat{d}^{q-1} + \widehat{d}^{q-2} d + \dots + \widehat{d} d^{q-2} + d^{q-1})\\
&\geq v_{\Gamma}(\widehat{d} - d) + (q-1) v_{\Gamma}(\widehat{d}) > q v_{\Gamma}(\widehat{d}) = v_{\Gamma}(\widehat{d}^q), \quad \mbox{hence $v_{\Gamma}(\tfrac{\widehat{d}^q - d^q}{\widehat{d}^q}) > 0$.}
\end{align*}

Consequently $v(\tfrac{\widehat{d}^q - d^q}{\widehat{d}^q}) > 0$, i.e.\ $v(\widehat{d}^q - d^q) > v(\widehat{d}^q)$; equivalently $v(c - d^q) > v(c)$ and hence $w(c - d^q) > w(c)$ as desired. \bs

\begin{thm}\label{mainthmforhens}
Let $q > \charac(k)$ be prime, $v$ a henselian valuation on $k((\Gamma))$ which factors through $v_{\Gamma}$, where $k$ and $\Gamma$ are countable, and $K_q$ as above. Then:
\begin{enumerate}
    \item $(K_q, v_{\Gamma}|_{K_q})$ is an immediate extension of $(k(\Gamma), v_{\Gamma})$;
    \item $(K_q, v|_{K_q})$ is an immediate extension of $(k(\Gamma), v)$;
    \item $\Z$ $($resp.\! $\mathbb{F}_p[z])$ is definable in $K_q$ in the language $\{0, 1, +, \times, \mathcal{O}_{v|_{K_q}}\}$;
    \item if $a \in k((\Gamma))\setminus K_q$ is \emph{pliant} over $K_q$, then $q | [K_q(a) : K_q]$.
\end{enumerate}
\end{thm}

\proof
\textit{(1) \& (2)} follow from the fact that $k(\Gamma) \subseteq K_q \subseteq k((\Gamma))$. For \textit{(3)}, $\Z$ (resp.\! $\mathbb{F}_p[z]$) is $\mathcal{L}_{val}$-{definable} in $K_q$ as $K_q \cap (k((\Gamma)))^q$ is $\mathcal{L}_{val}$-definable in $K_q$ by \emph{Lemma \ref{littledeflem}}, and this is sufficient to define $\Z$ (resp.\! $\mathbb{F}_p[z]$) by \emph{Lemma \ref{thelittlelem}}. Finally for \textit{(4)}, note for some $n$ sufficiently large, we have $a = a_{n}$ and $[E_{4n}(a_n) : E_{4n}] = [K_q(a) : K_q]$, as we assume $a \not\in K_q$. By construction (\textsc{Case 1}), $q|[E_{4n}(a) : E_{4n}]$ as desired. \bs

\begin{remark}\label{ultra}
(\cite[pp.\ 173--178]{EP}.) Let $S$ be an infinite set of indices and $\mathcal{U}$ a nonprincipal ultrafilter on $S$. For $s \in S$, let $(K_s, v_s)$ be a valued field. One may take an ultraproduct $\prod_{s \in S} (K_s, v_s)/\mathcal{U}$ of valued fields, and obtain a (valued) field $\mathbb{K} = \prod_{s \in S} K_s/\mathcal{U}$ with value group $\prod_{s \in S} v_s K_s/\mathcal{U}$ and residue field $\prod_{s \in S} K_s v_s/\mathcal{U}$, under the valuation $\prod v_s$ defined by:
\begin{gather*}
\textstyle\prod v_s ([(a_s)_{s \in S}]) = [(v_s(a_s))_{s \in S}], \quad \mbox{with residue}\\
res_{\textstyle\prod v_s} : \mathcal{O}_{\textstyle\prod v_s} \rightarrow \prod_{s \in S} K_s v_s/\mathcal{U};\quad [(x_s)] \mapsto [(res_{v_s}(x_s))].
\end{gather*}
(Here $a_s \in K_s$ for $s \in S$, and $[\cdot]$ represents the equivalence class of tuples modulo $\mathcal{U}$.) \sq
\end{remark}

\begin{cor}\label{thisis}
Let $(K, v)$ be an equicharacteristic $0$ henselian nontrivially valued field. $\TTh(K; \mathcal{L}_{val})$ is finitely undecidable. 

Moreover, if $\mathcal{O}_v$ is $\mathcal{L}_r$-definable, then $K$ is finitely undecidable as a field.
\end{cor}

\proof
Writing $k = Kv$ and $\Gamma = vK$, by \emph{Theorem \ref{akeref}} we have $(K, v) \equiv (k((\Gamma)), v_{\Gamma})$. By the \emph{Downwards L{\"o}wenheim-Skolem Theorem} we may also assume $k, \Gamma$ are countable. Set $v'$ (on $k$) to be the trivial valuation; in which case $v = v' \circ v_{\Gamma} = v_{\Gamma}$. By \emph{Lemma \ref{concrete}}, $(R, v_R)\! =\! (k(\Gamma), v_{\Gamma}|_{k(\Gamma)}),(Z, v_Z)\! = \!(k((\Gamma)), v_{\Gamma}),F\! =\! k((\Gamma)) \cap \widetilde{k(\Gamma)}$ satisfy \emph{Assumption $(\otimes)$}.

Let $q$ be prime and consider $(K_q, w = v_{\Gamma}|_{K_q}) \subseteq (k((\Gamma)), v_{\Gamma})$ arising from the \textsc{Modified Construction}, where now `pliant' means \emph{algebraic}. In particular $K_q$ is an equicharacteristic $0$ valued field with residue field $k$ and value group $\Gamma$. We will verify the henselianity axioms $\varphi_1, \dots, \varphi_{q-1}$ are satisfied, where $\varphi_n$ is
$$\forall a_0, \dots, a_{n-2} \left( \bigwedge_i a_i \in \mathfrak{m}_v \rightarrow \exists x\, [x^n + x^{n-1} + a_{n-2} x^{n-2} + \dots + a_0 = 0]\right).$$

Take $l \leq q-1$ and fix $a_0, \dots, a_{l-2} \in \mathfrak{m}_{w}$. Suppose $X^l + X^{l-1} + a_{l-2} X^{l-2} + \dots + a_0 = 0$ has no solution in $K_q$, and $\alpha \in k((\Gamma)) \setminus K_q$ satisfies this equation. Then $K_q(\alpha)/K_q$ is a finite proper extension. By \emph{Theorem \ref{mainthmforhens}}, $q| [K_q(\alpha) : K_q]$, however $q > l$ and $[K_q(\alpha) : K_q] \leq l$; a contradiction. We conclude that $K_q \models \varphi_l$ for all $l \leq q-1$; in particular for $n \geq 0$ fixed, $K_q \models \varphi_n$ for all primes $q > n$.

Let $\mathcal{U}$ be a nonprincipal ultrafilter on the set of primes, and let $\mathbb{K}$ be the ultraproduct $\prod_{q} K_q / \mathcal{U}$. By \emph{Remark \ref{ultra} \& \L o\'{s}' Theorem}, $\mathbb{K}$ is an equicharacteristic $0$ henselian valued field with residue field ($\mathcal{L}_r$-elementarily equivalent to) $k$, and value group ($\mathcal{L}_{oag}$-elementarily equivalent to) $\Gamma$. Hence by \emph{Theorem \ref{akeref}}, $\mathbb{K} \equiv_{\mathcal{L}_{val}} k((\Gamma))$. Thus $\mathbb{K} \models T$ for any finite subtheory $T \subseteq \Th(k((\Gamma)); \mathcal{L}_{val})$, and hence for some prime $q$, $K_q \models T$. {By \emph{Theorem \ref{mainthmforhens} \& Corollary \ref{appendixacor}}, $\Th(K_q; \mathcal{L}_{val})$ is hereditarily undecidable, making $T$ undecidable as required.} 

Finally, if $\mathcal{O}_v \subseteq K$ is $\mathcal{L}_r$-definable by the formula $\varpi(x, \overline{y})$ (where $\overline{y} = y_1, \dots, y_n$ denote parameter variables), then $\mathbb{K} \models \exists y_1, \dots, y_n \forall x\, (x \in \mathcal{O} \leftrightarrow \varpi(x, \overline{y}))$ too. For any finite subtheory $S \subseteq \Th(K; \mathcal{L}_r)$, by \emph{\L o\'{s}' Theorem} there exists a prime $l$ such that $K_l \models S \land \exists y_1, \dots, y_n \forall x\, (x \in \mathcal{O} \leftrightarrow \varpi(x, \overline{y}))$. By \emph{Theorem \ref{mainthmforhens} \& Corollary \ref{appendixacor}}, $\Th(K_l; \mathcal{L}_r)$ is hereditarily undecidable, making $S$ undecidable as desired. \bs

\bigskip
\section{Henselian Valued Fields: Further Discourse}\label{hvfagain}

We may extend the results of the previous section from equicharacteristic 0 henselian valued fields, to mixed characteristic henselian valued fields, using the \emph{standard decomposition}\footnote{This is the terminology used by Anscombe \& Jahnke \cite{anscombejahnke}; in the literature it is also known as the \emph{canonical decomposition}.}:

Let $(K, v)$ be a valued field of mixed characteristic $(0, p)$ with value group $\Gamma$. Define $\Delta_0$ to be the minimal convex subgroup of $\Gamma$ containing $v(p)$, and $\Delta_p$ to be the maximal convex subgroup of $\Gamma$ not containing $v(p)$. We will consider {the valuation(s) $v_0 : K \rightarrow \Gamma / \Delta_0$ (resp.\ $v_p : K \rightarrow \Gamma / \Delta_p$) corresponding to the coarsening(s) of $v$ with respect to $\Delta_0$ (resp.\ $\Delta_p$), and the induced valuation(s) $\hat{v}_0 : K v_0 \rightarrow \Delta_0$ (resp.\ $\hat{v}_p : K v_p \rightarrow \Delta_p$) with residue field(s) $Kv$. Also consider $\overline{v_p} : K v_0 \rightarrow \Delta_0 / \Delta_p$ with residue field $K v_p$;} this arises as the coarsening of $v_0$ with respect to $\Delta_p$ (as $\Delta_p < \Delta_0 \leq \Gamma$). These fit together in the following way:
\[
\begin{tikzcd}
K \arrow[r, "res_{v_0}"]  \arrow[d, "v_0"] & K v_0 \arrow[r, "res_{\overline{v_p}}"] \arrow[d, "\overline{v_p}"] & K v_p \arrow[d, "\hat{v}_p"] \arrow[r, "res_{\hat{v}_p}"] & Kv\\
\Gamma / \Delta_0 & \Delta_0 / \Delta_p & \Delta_p &
\end{tikzcd}
\]
\[
\begin{tikzcd}
K \arrow[r, "res_{v_0}"]  \arrow[d, "v_0"] & K v_0 \arrow[d, "\hat{v}_0"] \arrow[r, "res_{\hat{v}_0}"]& Kv && K \arrow[r, "res_{v_p}"]  \arrow[d, "v_p"] & K v_p \arrow[d, "\hat{v}_p"] \arrow[r, "res_{\hat{v}_p}"] & Kv\\
\Gamma/\Delta_0 &\Delta_0 &&& \Gamma / \Delta_p & \Delta_p &\\
\end{tikzcd}
\]

If $\Delta_0 = vK$ then by the \emph{Compactness Theorem} there is an elementary extension $(K, v) \prec (K^{\flat}, v^{\flat})$ containing an element $c^{\flat}$ such that $n \cdot v^{\flat}(p) < v^{\flat}(c^{\flat}) < \infty$ for all $n \geq 0$. Hence the minimal convex subgroup of $v^{\flat}K^{\flat}$ containing $v^{\flat}(p)$ does not contain $v^{\flat}(c^{\flat})$, i.e.\ $\Delta_0 < v^{\flat}K^{\flat}$. As $(K, v) \equiv (K^{\flat}, v^{\flat})$, for the purposes of proving finite undecidability we may assume WLOG $\Delta_0 < vK$. Consider:

\begin{lem}\label{aj65}
\emph{\cite[Lemma 6.5]{anscombejahnke}.} Let $T$ be a theory of bivalued fields $(K, v', v)$ with $v'$ an equicharacteristic $0$ henselian coarsening of $v$, and suppose that $T$ entails complete theories of valued fields $(K, v')$ and $(Kv', \overline{v})$. Then $T$ is complete. \bsn
\end{lem}

\begin{cor}\label{cor2}
Let $(K, v)$ be a mixed characteristic henselian nontrivially valued field, and $v_0$ as above. There exists a nontrivial ordered abelian group $\Omega$ such that $(K, v_0, v) \equiv (Kv_0((\Omega)), v_{\Omega}, \widehat{v}_0 \circ v_{\Omega})$.
\end{cor}

\proof Define $\Omega = vK/\Delta_0$. By \emph{Theorem \ref{akeref}}, $(K, v_0) \equiv (Kv_0((\Omega)), v_{\Omega})$. Since $v = \widehat{v}_0 \circ v_0$ in the standard decomposition of $(K, v)$,
\begin{align*}
    &(Kv_0, \widehat{v}_0) \qquad \mbox{(residue field of $(K, v_0)$ with the induced valuation)}\\
    =\,\, &(Kv_0, \widehat{v}_0) \qquad \mbox{(residue field of $(Kv_0((\Omega)), v_{\Omega})$ with the induced valuation).}
\end{align*}
By \emph{Lemma \ref{aj65}}, we conclude $(K, v_0, v) \equiv (Kv_0((\Omega)), v_{\Omega}, \widehat{v}_0 \circ v_{\Omega})$ as desired. \bs

\begin{thm}\label{newmixed}
Let $(K, v)$ be a mixed characteristic henselian nontrivially valued field. $\TTh(K; \mathcal{L}_{val})$ is finitely undecidable.

Moreover, if $\mathcal{O}_v$ is $\mathcal{L}_r$-definable, then $K$ is finitely undecidable as a field.
\end{thm}

\proof
As finite undecidability is a property of $\Th(K; \mathcal{L}_{val})$, we may assume WLOG $(K, v_0, v) = (Kv_0((\Omega)), v_{\Omega}, \widehat{v}_0 \circ v_{\Omega})$ where $Kv_0$ and $\Omega$ are countable, by the \emph{Downwards L{\"o}wenheim-Skolem Theorem \& Corollary \ref{cor2}}. Let $(Z, v_Z) = (Kv_0((\Omega)), \widehat{v}_0 \circ v_{\Omega})$, $(R, v_R) = (Kv_0(\Omega),\, \widehat{v}_0 \circ v_{\Omega})$ and $F = Z \cap \widetilde{R}$. By \emph{Lemma \ref{concrete}}, \emph{Assumption $(\otimes)$} is satisfied: for $q > p$ prime, denote by $K_q$ be the field given by the \textsc{Modified Construction} -- where now `pliant' means \emph{algebraic} -- with valuation $(\widehat{v}_0 \circ v_{\Omega})|_{K_q} = \widehat{v}_0 \circ v_{\Omega}|_{K_q}$. We have the following diagram of fields:
\[
\begin{tikzcd}
Kv_0((\Omega)) \arrow[r, "res_{v_{\Omega}}"] & Kv_0 \arrow[r, "res_{\widehat{v}_0}"] & Kv\\
K_q \arrow[u, dash]\arrow[r, "res_{v_{\Omega}}"] & Kv_0 \arrow[r, "res_{\widehat{v}_0}"] \arrow[u, equal] & Kv \arrow[u, equal]\\
Kv_0(\Omega) \arrow[u, dash] \arrow[r, "res_{v_{\Omega}}"] & Kv_0 \arrow[r, "res_{\widehat{v}_0}"] \arrow[u, equal] & Kv \arrow[u, equal]
\end{tikzcd}
\]
Consider $K_q$ as a multisorted structure:
$$\mathcal{K}_q = (K_q,\, Kv_0,\, Kv,\, \Omega,\, \Delta_0,\, vK;\,\, res_{v_{\Omega}|_{K_q}},\, res_{\widehat{v}_0},\, res_{(\widehat{v}_0 \circ v_{\Omega})|_{K_q}};\, v_{\Omega}|_{K_q},\, \widehat{v}_0,\, (\widehat{v}_0 \circ v_{\Omega})|_{K_q}),$$
which encompasses the diagram
\[
\begin{tikzcd}[column sep = large]
K_q \arrow[r, "res_{v_{\Omega}|_{K_q}}"] \arrow[dd, bend right, swap, "(\widehat{v}_0 \circ v_{\Omega})|_{K_q}"] \arrow[d, "v_{\Omega}|_{K_q}"] \arrow[rr, controls={+(80:1.25) and +(100:1.25)}, "res_{(\widehat{v}_0 \circ v_{\Omega})|_{K_q}}"] & Kv_0 \arrow[d, "\widehat{v}_0"]\arrow[r, "res_{\widehat{v}_0}"] & Kv\\
\Omega & \Delta_0 &\\
vK & &
\end{tikzcd}
\]

Let $\mathcal{U}$ be a nonprincipal ultrafilter on the set of primes larger than $p$, and let $\mathbb{K}$ be the ultraproduct $\prod_{q>p} \mathcal{K}_q/\mathcal{U}$. (We abuse notation to also denote the home sort of $\mathbb{K}$ by $\mathbb{K}$.) We have the following properties of $\mathbb{K}$:
\begin{itemize}
    \item $\mathbb{K}$ has valuation $\mathfrak{v}_0=\prod v_{\Omega}|_{K_q}$ with residue field $\prod_{q>p} Kv_0/\mathcal{U} = Kv_0^{\mathcal{U}}$ and value group $\prod_{q>p} \Omega/\mathcal{U} = \Omega^{\mathcal{U}}$. (This is \emph{Remark \ref{ultra}}.) Furthermore, $(\mathbb{K}, \mathfrak{v}_0)$ is a \textit{henselian} valued field. Indeed, we claim for $q > p$ prime the henselianity axioms $\varphi_1, \dots, \varphi_{q-1}$ are satisfied in $K_q$:
    
    Let $l < q$ and fix $a_0, \dots a_{l-2} \in \mathfrak{m}_{v_{\Omega}|_{K_q}}$. Suppose $X^l + X^{l-1} + a_{l-2}X^{l-2} + \dots + a_0 = 0$ has no solution in $K_q$ {-- though there exists a solution $\alpha \in K_q^h$, and $K_q^h \subseteq Kv_0((\Omega))$ by the universal property of henselisations \cite[Theorem 5.2.2]{EP}. Then $K_q(\alpha)/K_q$ is a finite proper extension.} By \emph{Theorem \ref{mainthmforhens}}, $q|[K_q(\alpha) : K_q]$, however $q > l$ and $[K_q(\alpha) : K_q] \leq l$; a contradiction.
    
    Therefore, by \textit{\L o\'{s}' Theorem} (i.e.\ \cite[Theorem 4.1.9]{changkeisler}) $(\mathbb{K}, \mathfrak{v}_0)$ is an equicharacteristic 0 henselian valued field with $\mathfrak{v}_0 \mathbb{K} \equiv_{\mathcal{L}_{oag}} \Omega$ and $\mathbb{K}\mathfrak{v}_0 \equiv_{\mathcal{L}_r} Kv_0$. By \emph{Theorem \ref{akeref}} we conclude $(\mathbb{K}, \mathfrak{v}_0) \equiv (Kv_0((\Omega)), v_{\Omega})$.
    \item $Kv_0^{\mathcal{U}}$ has valuation $\widehat{v}_0^{\mathcal{U}} = \prod \widehat{v}_0$ with residue field $vK^{\mathcal{U}}$ and value group $\Delta_0^{\mathcal{U}}$ (\emph{Remark \ref{ultra}}). By \cite[Theorem 4.1.9]{changkeisler}, $(Kv_0^{\mathcal{U}}, \widehat{v}_0^{\mathcal{U}}) \equiv (Kv_0, \widehat{v}_0)$.
    \item Hence $\mathbb{K}$ can be equipped with a valuation $\mathfrak{v}_1 := \widehat{v}_0^{\mathcal{U}} \circ \mathfrak{v}_0$, and $\mathfrak{v}_0$ is an equicharacteristic $0$ henselian coarsening of $\mathfrak{v}_1$. $\mathbb{K}$ also has a valuation $\mathfrak{v}_2=\prod (\widehat{v}_0 \circ v_{\Omega})|_{K_q}$, and we claim $\mathcal{O}_{\mathfrak{v}_1} = \mathcal{O}_{\mathfrak{v}_2}$, so $(\mathbb{K}, \mathfrak{v}_1) \equiv (\mathbb{K}, \mathfrak{v}_2)$. Indeed,
    \begin{align*}
        x = [(x_q)] \in \mathcal{O}_{\mathfrak{v}_1} &\iff [(x_q)] \in res_{\mathfrak{v}_0}^{-1}(\mathcal{O}_{\widehat{v}_0^{\mathcal{U}}})\\
        &\iff res_{\mathfrak{v}_0}([(x_q)]) = [(res_{v_{\Omega}|_{K_q}}(x_q))] \in \mathcal{O}_{\widehat{v}_0^{\mathcal{U}}}\\
        &\iff \{q \mbox{ : } res_{v_{\Omega}|_{K_q}}(x_q) \in \mathcal{O}_{\widehat{v}_0}\} \in \mathcal{U}\\
        &\iff \{q \mbox{ : } x_q \in \mathcal{O}_{(\widehat{v}_0 \circ v_{\Omega})|_{K_q}}\} \in \mathcal{U} \iff [(x_q)] \in \mathcal{O}_{\mathfrak{v}_2}.
    \end{align*}
\end{itemize}
Considering $(\mathbb{K}, \mathfrak{v}_0, \mathfrak{v}_1)$ as a bivalued field, by \emph{Lemma \ref{aj65}} 
$$(\mathbb{K}, \mathfrak{v}_0, \mathfrak{v}_1) \equiv (Kv_0((\Omega)), v_{\Omega}, \widehat{v}_0 \circ v_{\Omega}) = (K, v_0, v).$$

Taking a reduct of the language, $(\mathbb{K}, \mathfrak{v}_2) \equiv (\mathbb{K}, \mathfrak{v}_1) \equiv (K, v)$. Therefore if $T$ is a finite subtheory of $\mbox{\ttfamily Th}(K; \mathcal{L}_{val})$, $(\mathbb{K}, \mathfrak{v}_2) \models T$, and hence for some $q$ sufficiently large, $(K_q, \widehat{v}_0 \circ v_{\Omega}|_{K_q}) \models T$. {By \emph{Theorem \ref{mainthmforhens} \& Corollary \ref{appendixacor}}, $\Th(K_q; \mathcal{L}_{val})$ is hereditarily undecidable, making $T$ undecidable as required.}

Finally, if $\mathcal{O}_v \subseteq K$ is $\mathcal{L}_r$-definable by the formula $\varpi(x, \overline{y})$ (where $\overline{y} = y_1, \dots, y_n$ denote parameter variables), then $\mathbb{K} \models \exists y_1, \dots, y_n \forall x\, (x \in \mathcal{O}_{\mathfrak{v}_2} \leftrightarrow \varpi(x, \overline{y}))$ too. For any finite subtheory $S \subseteq \Th(K; \mathcal{L}_r)$, by \emph{\L o\'{s}' Theorem} there exists a prime $l$ such that $K_l \models S \land \exists y_1, \dots, y_n \forall x\, (x \in \mathcal{O} \leftrightarrow \varpi(x, \overline{y}))$. By \emph{Theorem \ref{mainthmforhens} \& Corollary \ref{appendixacor}}, $\Th(K_l; \mathcal{L}_r)$ is hereditarily undecidable, making $S$ undecidable as desired. \bs

What remains is to handle the case of equicharacteristic $p > 0$ henselian nontrivially valued fields. We will show this gap can be eliminated for certain fields of this type (e.g.\ those fields also satisfying \emph{NIP}, by using a nice algebraic classification of such fields by Anscombe \& Jahnke \cite[Theorem 5.1]{anscombejahnke}). We require the following definitions:

\begin{definition}
A valued field $(K, v)$ is said to be \emph{(separably) defectless} if whenever $L/K$ is a finite (separable) extension, $[L : K] = \sum_{w \supseteq v} e(w/v) f(w/v)$, where $w$ ranges over all prolongations of $v$ to $L$, $e(w/v) = (wL : vK)$ is the ramification degree and $f(w/v) = [Lw : Kv]$ is the inertia degree of the valued field extension $(L,w)/(K,v)$. 
\end{definition} 

\begin{definition}\label{kaplansky}
A valued field $(K, v)$ of residue characteristic $p>0$ is \textit{Kaplansky}\footnote{Equivalently (\cite[Remark 2.2]{anscombejahnke}) $(K, v)$ of residue characteristic $p>0$ is Kaplansky if and only if $vK$ is $p$-divisible and $Kv$ admits no \textit{finite} proper extensions of degree divisible by $p$.} if $vK$ is $p$-divisible, and $Kv$ is perfect and admits no proper separable algebraic extensions of degree divisible by $p$. 
\end{definition}

\begin{thm}\label{sylvyfranzi}
\emph{\textbf{(Anscombe-Jahnke)}} Let $(K, v)$ be a positive equicharacteristic NIP henselian nontrivially valued field. Then $(K, v)$ is separably defectless Kaplansky.
\end{thm}

\proof
An immediate consequence of \cite[Proposition 3.1]{anscombejahnke}. \bs

{A corollary to this, noted by Anscombe \& Jahnke, is that by Delon \cite[Th{\'e}or{\`e}me 3.1]{delon}} the $\mathcal{L}_{val}$-theory of equicharacteristic $p > 0$ henselian separably defectless valued fields $(K,v)$ of imperfection degree $e$, with residue field $\mathcal{L}_r$-elementarily equivalent to $Kv$ and value group $\mathcal{L}_{oag}$-elementarily equivalent to $vK$, is complete. (As we are concerned with \emph{finitely axiomatised} subsets of {$\Th(K; \mathcal{L}_{val})$, we will assume WLOG $e < \infty$.}) 

\bigskip
{For a valued field $(B, v_B)$ to be separably defectless, it is sufficient for it to satisfy the first-order $\mathcal{L}_{val}$-statements\footnote{The argument of \cite[Lemma 2.4]{anscombejahnke} confirms $(\raisebox{-0.5mm}{$\vardiamond$}_M)$ is a first-order $\mathcal{L}_{val}$-statement (with ``defectless'' replaced by ``separably defectless'', and field extensions made separable where appropriate).} $(\raisebox{-0.5mm}{$\vardiamond$}_M)$ for all $M \geq 1$:

\begin{align*}
  &\mbox{For all finite separable extensions $D/B$ of degree $\leq M$, the equality}\\
  &\hspace{37mm}[D : B] = \sum_{w \supseteq v_B} e(w/v_B) f(w/v_B)\hspace{7mm} \tag{$\raisebox{-0.5mm}{$\vardiamond$}_M$}\\
  &\mbox{holds, where $w$ ranges over all prolongations of $v_B$ to $D$, $e(w/v_B)$ is}\\
  &\mbox{the ramification index and $f(w/v_B)$ is the inertia degree.}
\end{align*}

\vspace{2mm}
\begin{lem}\label{defectclaim}
Fix $M \in \N_{>0}$ and $q > M^M$ prime. Assume $(R, v_R)$, $(Z, v_Z)$, and $F$ satisfy \emph{Assumption $(\otimes)$} and let $K_q$ be the field resulting from the \textsc{Modified Construction}. If $Z$ is separably defectless Kaplansky, then $K_q \models (\raisebox{-0.5mm}{$\vardiamond$}_M)$.
\end{lem}}

\proof
Let $D/K_q$ be a separable extension of degree $\leq M$. By taking the normal closure may assume $D/K_q$ is Galois and of degree $d\leq M^M$. If $D/K_q$ is proper, by degree reasons $Z \cap D = K_q$. Indeed, if $Z \supset D' \supset K_q$ is a separable extension of degree $\leq M^M$, by the \emph{Primitive Element Theorem} there exists $\alpha \in Z \cap K_q^s$ such that $D' = K_q(\alpha)$. By \emph{Theorem \ref{mainthmforhens}}, if $\alpha \not\in K_q$ then $q | [K_q(\alpha) : K_q]$, however $q > M^M$ and $[K_q(\alpha) : K_q] \leq M^M$. This is a contradiction, hence $\alpha \in K_q$. Thus $Z$ and $D$ are linearly disjoint over $K_q$. Consider the following tower of extensions: 
\[
\begin{tikzcd}[column sep=normal, row sep=tiny]
& Z D\arrow[dd, dash] \\
Z \arrow[dd, dash] \arrow[ur, dash, "d", pos=0.6]  \\
& D   \\
K_q  \arrow[ur, dash, "d", pos=0.6]
\end{tikzcd}
\]

By the tower property (\cite[Lemma 2.5.3]{friedjarden}), $D$ is linearly disjoint to $K_q^h$, a subfield of $Z$ by the universal property of henselisations. Thus $w = v_Z|_{K_q}$ extends uniquely to $(D, w')$, by e.g.\ \cite[Lemma 2.1]{abfvk}. As $(Z, v_Z)$ is henselian, $v_Z$ extends uniquely to $(ZD, v_Z')$, and restricts to $(D, w')$. {As $Z/K_q$ is immediate, as the valuation extensions are unique, and as $Z$ is linearly disjoint from $D$ over $K_q$, 
\begin{align*}
p^{\nu} e(w'/w) f(w'/w) &= [D : K_q] \qquad\mbox{where $\nu \geq 0$, by \cite[Theorem 3.3.3]{EP},}\\
&= [ZD : Z] \qquad \mbox{by \cite[Corollary 2.5.2]{friedjarden},}\\
&= e(v_Z'/v_Z) f(v_Z' / v_Z) \qquad\mbox{as $Z$ is separably defectless}.
\end{align*}
In addition, as $Z$ is Kaplansky, $p \nmid e(v_Z'/v_Z), f(v_Z' / v_Z)$. Hence 
$$[D : K_q] = e(w'/w) f(w'/w) = \sum_{w' \supseteq w} e(w'/w) f(w'/w).$$
We conclude} $K_q \models (\raisebox{-0.5mm}{$\vardiamond$}_M)$ as desired. \bs

\noindent We are ready to prove:

\begin{thm}\label{nipequichar}
Let $(K, v)$ be an equicharacteristic $p > 0$ separably defectless Kaplansky henselian nontrivially valued field. $\TTh(K; \mathcal{L}_{val})$ is finitely undecidable.

Moreover, if $\mathcal{O}_v$ is $\mathcal{L}_r$-definable, then $K$ is finitely undecidable as a field.
\end{thm}

\proof
{Fix $e \in \N$, the imperfection degree of $K$.} Let $k$, $\Gamma$ be countable models of $\Th(Kv; \mathcal{L}_r)$, $\Th(vK; \mathcal{L}_{oag})$. By definition $k$ is perfect and $\Gamma$ is $p$-divisible, and thus $k(\Gamma)$ and $k((\Gamma))$ are perfect. {By \emph{Lemma \ref{concrete}} we may choose $t_1, \dots, t_e \in k((\Gamma))$ transcendental and algebraically independent over $k(\Gamma)$} and consider the field $k(\Gamma)(t_1, \dots, t_{e-1})$ -- by \cite[Lemma 2.7.2]{friedjarden} its imperfection degree is exactly $e-1$. Finally, set $F = k((\Gamma)) \cap (k(\Gamma)(t_1, \dots, t_{e-1}))^s$. With $(R, v_R) = (k(\Gamma), v_{\Gamma})$, $(Z, v_Z) = (k((\Gamma)), v_{\Gamma})$, \emph{Assumption $(\otimes)$} is satisfied by \emph{Lemma \ref{concrete}}. Note the imperfection degree of $F$ is exactly $e-1$, and $(Z, v_Z)$ is perfect separably defectless Kaplansky.

As there exists an element $t \in k((\Gamma))$ transcendental over $F$, by \emph{Theorem \ref{mainthmforhens}} there exists a field $F(t) \subseteq K_q \subseteq k((\Gamma))$ such that if $w = v_{\Gamma}|_{K_q}$, then $K_q w = k$, $wK_q = \Gamma$, {$\mathbb{F}_p[z]$ is $\mathcal{L}_{val}$-definable} in $K_q$ where $z \in k(\Gamma)$ is transcendental over $\mathbb{F}_p$, and if $a \in k((\Gamma))\setminus K_q$ is \textit{pliant} over $K_q$, then $q | [K_q(a) : K_q]$. 
\begin{itemize}[leftmargin = 2cm]
\item[If $e > 0$:] `pliant' will refer to \emph{separably algebraic} elements. By the construction of $K_q$ as a union of separable extensions of $F(t)$, $K_q/k(\Gamma)(t_1, \dots, t_{e-1}, t)$ is separably algebraic and has degree of imperfection $e$ by \cite[Lemma 2.7.3]{friedjarden}.
\item[If $e = 0$:] `pliant' will refer to \emph{algebraic} elements. As noted in \emph{Remark \ref{littleprob}}, $K_q$ is perfect.
\end{itemize}

In either case, {\emph{$K_q$ is an equicharacteristic $p$ nontrivially valued field of imperfection degree $e$, with residue field $k$ and value group $\Gamma$.}} We will verify the henselianity axioms $\varphi_1, \dots, \varphi_{q-1}$ are satisfied. Let $l < q$ and fix $a_0, \dots a_{l-2} \in \mathfrak{m}_{w}$. Suppose $X^l + X^{l-1} + a_{l-2}X^{l-2} + \dots + a_0 = 0$ has no solution in $K_q$ -- however there exists a solution $\alpha \in K_q^h$, and $K_q^h \subseteq k((\Gamma))$ by the universal property of henselisations \cite[Theorem 5.2.2]{EP}. Then $K_q(\alpha)/K_q$ is a finite proper extension; moreover, $\alpha$ is pliant over $K_q$ as $K_q(\alpha) \subseteq K_q^h \subseteq K_q^s \subseteq \widetilde{K_q}$. By \emph{Theorem \ref{mainthmforhens}}, $q|[K_q(\alpha) : K_q]$, however $q > l$ and $[K_q(\alpha) : K_q] \leq l$; a contradiction.

{Let $Q$ be the set of primes $q > p$ and $\mathcal{U}$ a nonprincipal ultrafilter on $Q$. Let $\mathbb{K} = \prod_{q \in Q} K_q / \mathcal{U}$; by \textit{\L o\'{s}' Theorem}, $\mathbb{K}$ is an equicharacteristic $p$ henselian nontrivially valued field, of imperfection degree $e$, and residue field $\mathcal{L}_r$-elementarily equivalent to $Kv$ and value group $\mathcal{L}_{oag}$-elementarily equivalent to $vK$. Furthermore, $\mathbb{K}$ is separably defectless. Indeed, given $M \in \mathbb{N}_{>0}$, $K_q \models (\raisebox{-0.5mm}{$\vardiamond$}_M)$ for all primes $q > M^M$ by \emph{Lemma \ref{defectclaim}}. Hence, by \textit{\L o\'{s}' Theorem}, $\mathbb{K} \models (\raisebox{-0.5mm}{$\vardiamond$}_M)$ for all $M \geq 1$. By \cite[Th{\'e}or{\`e}me 3.1]{delon}, $\mathbb{K} \equiv_{\mathcal{L}_{val}}\! K$.} 

Let $T$ be a finite subtheory of $\Th(K; \mathcal{L}_{val})$. As $\mathbb{K} \models T$, for some $q \in Q$ we have $K_q \models T$. By \emph{Theorem \ref{mainthmforhens} \& Corollary \ref{appendixacor}}, {$\Th(K_q; \mathcal{L}_{val})$ is hereditarily undecidable, hence $T$ is undecidable as required.} 

Finally, if $\mathcal{O}_v \subseteq K$ is $\mathcal{L}_r$-definable by the formula $\varpi(x, \overline{y})$ (where $\overline{y} = y_1, \dots, y_n$ denote parameter variables), then $\mathbb{K} \models \exists y_1, \dots, y_n \forall x\, (x \in \mathcal{O} \leftrightarrow \varpi(x, \overline{y}))$ too. For any finite subtheory $S \subseteq \Th(K; \mathcal{L}_r)$, by \emph{\L o\'{s}' Theorem} there exists a prime $l$ such that $K_l \models S \land \exists y_1, \dots, y_n \forall x\, (x \in \mathcal{O} \leftrightarrow \varpi(x, \overline{y}))$. By \emph{Theorem \ref{mainthmforhens} \& Corollary \ref{appendixacor}}, $\Th(K_l; \mathcal{L}_r)$ is hereditarily undecidable, making $S$ undecidable as desired. \bs



\begin{cor}\label{finalcor}
\phantom{boo}
\begin{enumerate}
    \item[\textit{(1)}] If $(K, v)$ is an NIP henselian nontrivially valued field, $\TTh(K; \mathcal{L}_{val})$ is finitely undecidable. Furthermore if $\mathcal{O}_v$ is $\mathcal{L}_r$-definable in $K$, then $K$ is finitely undecidable as a field.
    \item[\textit{(2)}] Every infinite dp-finite field is finitely undecidable.
    \item[\textit{(3)}] Assuming the \emph{NIP Fields Conjecture}, {every} infinite NIP field is finitely undecidable.
\end{enumerate}
\end{cor}

\proof
\phantom{boo}
\begin{enumerate}
    \item If $(K, v)$ is equicharacteristic 0 or mixed characteristic, this is a result of \emph{Corollary \ref{thisis}/Theorem \ref{newmixed}}. If $(K, v)$ is positive equicharacteristic then by Anscombe-Jahnke it is separably defectless Kaplansky (\emph{Theorem \ref{sylvyfranzi}}) hence its $\mathcal{L}_{val}$-theory is finitely undecidable by \emph{Theorem \ref{nipequichar}}.
    \item By the work of Johnson\footnote{This paper is the conclusion of the series \cite{johnsonia, johnsonib, johnsonii, johnsoniii,johnsoniv, johnsonv} by Johnson on dp-finite fields.} \cite[Corollary 4.16]{johnsonvi}, every infinite dp-finite field $K$ is either algebraically closed, real closed, or admits a non-trivial definable henselian valuation. If $K$ is algebraically or real closed, it is finitely undecidable by \cite[Corollary, p.\! 270]{ziegler}. Otherwise $K$ has an $\mathcal{L}_r$-definable valuation $v$. As $K$ is dp-finite, $\Th(K; \mathcal{L}_r)$ is NIP (following from \cite[Observation 4.13]{simon}), and hence $\Th(K; \mathcal{L}_{val})$ is NIP as the valuation is $\mathcal{L}_r$-definable. Thus $K$ is finitely undecidable, from (1).
    \item Recall \emph{Theorem \ref{maybejochen}}, where assuming the \textit{NIP Fields Conjecture} one can conclude every infinite NIP field $K$ is either separably closed (hence finitely undecidable by \cite[Corollary, p.\! 270]{ziegler}/\emph{Corollaries \ref{scf1} \&\ \ref{scf2}}), real closed (hence finitely undecidable by \cite[Corollary, p.\! 270]{ziegler}), or admits a nontrivial $\emptyset$-$\mathcal{L}_r$-definable henselian valuation. In this case, $K$ is finitely undecidable by (1).\bs
\end{enumerate}

\begin{exmp}\label{finalex}
We present some finitely undecidable fields, and some open questions, as a consequence of this work. 
\begin{enumerate}
    \item Let $(K, v)$ be any algebraic valued field extension of $(\Q_p, v_p)$ with non-divisible value group. By \cite[Lemma 3.6]{jochengalgroup}, the valuation is $\emptyset$-$\mathcal{L}_r$-definable. By \emph{Theorem \ref{newmixed}}, $K$ is a finitely undecidable field.
    \item It is unknown if $\mathbb{F}_p((t))$ is finitely undecidable. (Equally one can consider the finite undecidability of $\Th(\mathbb{F}_p((t)); \mathcal{L}_{val})$, as the valuation ring $\mathbb{F}_p[[t]]$ is $\emptyset$-$\mathcal{L}_r$-definable in $\mathbb{F}_p((t))$ by e.g.\ \cite[Lemma 3.6]{jochengalgroup}.) More generally if $F\!/\mathbb{F}_p$ is any algebraic extension, the finite undecidability of $F((t))$ -- or of $\Th(F((t)); \mathcal{L}_{val})$ -- is an open question.
\end{enumerate}

Recall the notion of \emph{t-henselianity}: a field is \emph{t-henselian} if it can be equipped with a topology compatible with the field operations, behaving very much like the topology arising from a nontrivial valuation, satisfying a `topological' Hensel's Lemma. This notion\linebreak was introduced by Prestel \& Ziegler \cite{prestelziegler}, and the topology shown to be $\mathcal{L}_r$-definable in nonseparably closed fields in \cite[p.\! 203]{presteltop}. Using saturation we may equally state (as Anscombe \& Jahnke \cite[p.\! 872]{ansjahn} do) that \emph{a field is t-henselian if it is $\mathcal{L}_r$-elementarily equivalent to a field $L$ which admits a nontrivial henselian valuation}, i.e.\ has a subset $\mathcal{O} \subseteq L$ that satisfies the definition of a nontrivial henselian valuation ring. (Cf.\ the discussion in \cite[\S 7]{prestelziegler} \& \cite[\S 3.4]{jochengalgroup}.) In some cases we can use t-henselianity -- a purely {field-theoretic}, $\mathcal{L}_r$-elementary property -- to conclude finite undecidability of a field.
\begin{itemize}
    \item[(3)] Let $K$ be a characteristic 0 t-henselian field with \emph{nonuniversal}\footnote{A profinite group $G$ is \emph{universal} if every finite group occurs as the image of a continuous morphism from $G$.} absolute Galois group. From\footnote{Cf.\ \cite[Theorem 3.15]{jochenjahnke} and the \textit{Remark loc.\! cit.}} \cite[Theorem 3.2.3]{jahnkethesis}, $K$ is either real closed, separably (= algebraically) closed, or henselian with respect to a nontrivial $\emptyset$-$\mathcal{L}_r$-definable valuation. Thus $K$ is finitely undecidable, by \cite[Corollary, p.\! 270]{ziegler} or \emph{Corollary \ref{thisis}} or \emph{Theorem \ref{newmixed}}.

    \item[(4)] Let $K$ be a positive characteristic \textit{NIP} t-henselian field. If $K$ is separably closed, it is finitely undecidable by \emph{Corollary \ref{scf1}/\ref{scf2}}. Otherwise let $L$ be a field with nontrivial henselian valuation $v$ such that $L \equiv K$; note $L$ is not separably closed and NIP \emph{as a field}. By \cite[Corollary 3.18]{jochenjahnke}, $L$ admits a nontrivial $\emptyset$-$\mathcal{L}_r$-definable henselian valuation. Hence $K$ is finitely undecidable, by \emph{Corollary \ref{finalcor} (1)}.
\end{itemize}
Two concrete examples of fields from (3) \& (4) are $\C((t))$ \& $\widetilde{\mathbb{F}_p}((\tfrac{1}{p^{\infty}}\!\Z))$, respectively.
\begin{itemize}
    \item[(5)] Boissonneau \cite{blaise} has recently extended the Anscombe-Jahnke results towards \emph{$n$-dependent} valued fields. The class of $n$-dependent theories was introduced by Shelah in \cite[\S 5 (H)]{shelahndep} (see \cite{artemetal} for further discussion), and is a proper generalisation of NIP (which corresponds to \emph{1-dependent}). $n$-dependence in groups and fields has been studied extensively by Hempel and Chernikov-Hempel \cite{hempel,chernhemp}.

    If $(K, v)$ is a positive characteristic $n$-dependent nontrivially valued field, it is henselian (by \cite[Theorem 3.1]{chernhemp}) and separably defectless Kaplansky (by \cite[Lemma 3.10]{blaise}). Thus $\Th(K; \mathcal{L}_{val})$ is finitely undecidable, by \emph{Theorem \ref{nipequichar}}. Furthermore, by \cite[Proposition 8.4]{hempel}, $K$ is either separably closed (hence finitely undecidable, by \emph{Corollary \ref{scf1}/\ref{scf2}}) or admits an $\emptyset$-$\mathcal{L}_r$-definable nontrivial henselian valuation (separably defectless Kaplansky, by \cite[Lemma 3.10]{blaise}). We conclude from \emph{Theorem \ref{nipequichar}} if $K$ is a positive characteristic $n$-dependent t-henselian field, it is finitely undecidable. \sq
\end{itemize}
\end{exmp}

One may ask: \textit{is every nontrivially henselian valued field, which is finitely undecidable as an $\mathcal{L}_{val}$-structure, finitely undecidable} as a field\emph{?} This cannot be concluded easily from \emph{Corollary \ref{thisis}/Theorem \ref{newmixed}/Theorem \ref{nipequichar}}, as there exist henselian valued fields which do not admit \emph{any} nontrivial $\mathcal{L}_r$-definable henselian valuation, such as the Jahnke-Koenigsmann example \cite[Example 6.2]{jochenjahnkeagain}. Anscombe \& Jahnke discuss ``henselianity in the language of rings'' further in \cite{ansjahn}.


\bigskip
\section*{Acknowledgements}
The author extends his thanks to Professor Ehud Hrushovski, Professor Jochen\linebreak Koenigsmann, Professor Arno Fehm, Dr.\ Sylvy Anscombe and Dr.\ Konstantinos Kartas for many helpful conversations. Further thanks is due to S.\ Anscombe for her excellent suggestions on presentation, organisation, and abstraction.

\bigskip
\bibliography{refs}   

\begin{thebibliography}{10}

\bibitem{ansjahn}
{\sc Anscombe, S., and Jahnke, F.}
\newblock {Henselianity in the language of rings}.
\newblock {\em {Ann.\ Pure Appl.\ Logic} 169}, 9 (2018), 872--895.

\bibitem{anscombejahnke}
{\sc Anscombe, S., and Jahnke, F.}
\newblock {Characterising NIP henselian fields}.
\newblock {\em \emph{\href{https://arxiv.org/abs/1911.00309}{\ttfamily
  arXiv:1911.00309}}\/} (2019).

\bibitem{axkochen}
{\sc Ax, J., and Kochen, S.}
\newblock {Diophantine Problems Over Local Fields I}.
\newblock {\em {Amer.\ Jour.\ Math.} 87}, 3 (1965), 605--630.

\bibitem{abfvk}
{\sc {Blaszczok, A., and Kuhlmann, {F.-V.}}}
\newblock {On maximal immediate extensions of valued fields}.
\newblock {\em {Math. Nachr.} 290}, 1 (2017), 7--18.

\bibitem{blaise}
{\sc Boissonneau, B.}
\newblock Artin-schreier extensions and combinatorial complexity in henselian
  valued fields.
\newblock {\em \emph{\href{https://arxiv.org/abs/2108.12678v2}{\ttfamily
  arXiv:2108.12678}}\/} (2022).

\bibitem{changkeisler}
{\sc Chang, C., and Keisler, H.}
\newblock {\em {Model Theory}}.
\newblock Dover, 2012.
\newblock {Third Edition}.

\bibitem{cherlin-shelah}
{\sc Cherlin, G., and Shelah, S.}
\newblock {Superstable fields and groups}.
\newblock {\em Ann.\ Math.\ Logic 18}, 3 (1980), 227--270.

\bibitem{chernhemp}
{\sc Chernikov, A., and Hempel, N.}
\newblock {On n-dependent groups and fields II}.
\newblock {\em Forum Math. Sigma 9}, e38 (2021), 1--51.

\bibitem{artemetal}
{\sc Chernikov, A., Palacin, D., and Takeuchi, K.}
\newblock {On $n$-Dependence}.
\newblock {\em Notre Dame J. Form. Log. 60}, 2 (2019), 195--214.

\bibitem{delon}
{\sc Delon, F.}
\newblock {\em {\it Quelques propri{\'e}t{\'e}s des corps valu{\'e}s en
  th{\'e}orie des mod{\`e}les}}.
\newblock PhD thesis, {Universit{\'e} Paris VII}, 1981.

\bibitem{bouscaren}
{\sc Delon, F.}
\newblock {Separably closed fields}.
\newblock In {\em Model Theory and Algebraic Geometry}, E.~Bouscaren, Ed.
  Springer, 1999, pp.~143--177.
\newblock Lecture Notes in Mathematics, vol. 1696.

\bibitem{vdd}
{\sc van~den Dries, L.}
\newblock {Lectures on the Model Theory of Valued Fields}.
\newblock In {\em Model theory in Algebra, Analysis and Arithmetic}, H.~D.
  Macpherson and C.~Toffalori, Eds. Springer-Verlag, 2014, pp.~55--157.

\bibitem{EP}
{\sc Engler, A.~J., and Prestel, A.}
\newblock {\em {Valued Fields}}.
\newblock Springer, 2005.
\newblock {Springer Monographs in Mathematics}.

\bibitem{ershovake}
{\sc Ershov, Y.}
\newblock {On the elementary theory of maximal normed fields}.
\newblock {\em {Soviet Math.\ Dokl.} 6\/} (1965), 1390--1393.

\bibitem{eltt}
{\sc Ershov, Y., Lavrov, I., Taimanov, A., and Taitslin, M.}
\newblock {Elementary Theories}.
\newblock {\em Russian Math.\ Surveys 20}, 4 (1965), 35--105.
\newblock (English version).

\bibitem{friedjarden}
{\sc Fried, M., and Jarden, M.}
\newblock {\em {Field Arithmetic}}.
\newblock Springer-Verlag, 2008.
\newblock {Ergebnisse der Mathematik und ihrer Grenzgebiete. 3.\ Folge; A
  Series of Modern Surveys in Mathematics, Volume 11, 3rd Edition}.

\bibitem{hhj}
{\sc Halevi, Y., Hasson, A., and Jahnke, F.}
\newblock {Definable $V$-topologies, Henselianity and NIP}.
\newblock {\em J.\ Math.\ Log. 20}, 2 (2020), 1--33.

\bibitem{halevi-palacin}
{\sc Halevi, Y., and Palac\'{i}n, D.}
\newblock {The dp-rank of Abelian groups}.
\newblock {\em J.\ Symb.\ Logic 84}, 3 (2019), 957--986.

\bibitem{hempel}
{\sc Hempel, N.}
\newblock {On $n$-dependent groups and fields}.
\newblock {\em Math. Log. Q. 62}, 3 (2016), 215--224.

\bibitem{jahnkethesis}
{\sc Jahnke, F.}
\newblock {\em {\it Definable Henselian Valuations and Absolute Galois
  Groups}}.
\newblock PhD thesis, {University of Oxford}, 2014.

\bibitem{jochenjahnke}
{\sc Jahnke, F., and Koenigsmann, J.}
\newblock {Definable Henselian Valuations}.
\newblock {\em J.\ Symb.\ Logic 80}, 1 (2015), 85--99.

\bibitem{jochenjahnkeagain}
{\sc Jahnke, F., and Koenigsmann, J.}
\newblock {Defining Coarsenings of Valuations}.
\newblock {\em Proc. Edinb. Math. Soc. 60}, 3 (2017), 665--687.

\bibitem{johnsonii}
{\sc Johnson, W.}
\newblock {Dp-finite fields II: the canonical topology and its relation to
  henselianity}.
\newblock {\em \emph{\href{https://arxiv.org/abs/1910.05932}{\ttfamily
  arXiv:1910.05932}}\/} (2019).

\bibitem{johnsoniii}
{\sc Johnson, W.}
\newblock {Dp-finite fields III: inflators and directories}.
\newblock {\em \emph{\href{https://arxiv.org/abs/1911.04727}{\ttfamily
  arXiv:1911.04727}}\/} (2019).

\bibitem{johnsoniv}
{\sc Johnson, W.}
\newblock {Dp-finite fields IV: the rank 2 picture}.
\newblock {\em \emph{\href{https://arxiv.org/abs/2003.09130}{\ttfamily
  arXiv:2003.09130}}\/} (2020).

\bibitem{johnsonv}
{\sc Johnson, W.}
\newblock {Dp-finite fields V: topological fields of finite weight}.
\newblock {\em \emph{\href{https://arxiv.org/abs/2004.14732}{\ttfamily
  arXiv:2004.14732}}\/} (2020).

\bibitem{johnsonvi}
{\sc Johnson, W.}
\newblock {Dp-finite fields VI: the dp-finite Shelah conjecture}.
\newblock {\em \emph{\href{https://arxiv.org/abs/2005.13989}{\ttfamily
  arXiv:2005.13989}}\/} (2020).

\bibitem{johnsonia}
{\sc Johnson, W.}
\newblock {Dp-finite fields I(A): The infinitesimals}.
\newblock {\em Ann.\ Pure Appl.\ Logic 172}, 6 (2021), 102947.

\bibitem{johnsonib}
{\sc Johnson, W.}
\newblock {Dp-finite fields I(B): Positive characteristic}.
\newblock {\em Ann.\ Pure Appl.\ Logic 172}, 6 (2021), 102949.

\bibitem{jtwy}
{\sc Johnson, W., Tran, C.-M., Walsberg, E., and Ye, J.}
\newblock The \'{e}tale-open topology and the stable field conjecture.
\newblock {\em \emph{\href{https://arxiv.org/abs/2009.02319v2}{\ttfamily
  arXiv:2009.02319v2}}\/} (2021).

\bibitem{kaplan}
{\sc Kaplan, I., Onshuus, A., and Usvyatsov, A.}
\newblock Additivity of the dp-rank.
\newblock {\em Trans. Amer. Math. Soc. 365}, 11 (2013), 5783--5804.

\bibitem{jochengalgroup}
{\sc Koenigsmann, J.}
\newblock {Elementary characterization of fields by their absolute Galois
  group}.
\newblock {\em Siberian Adv.\ Math. 14}, 3 (2004), 16--42.

\bibitem{ober}
{\sc Koenigsmann, J., Pasten, H., Shlapentokh, A., and Vidaux, X.}
\newblock {Definability and Decidability Problems in Number Theory}.
\newblock {\em Oberwolfach Rep. 13}, 4 (2016), 2793--2866.

\bibitem{krupinsky-pillay}
{\sc Krupi\'{n}ski, K., and Pillay, A.}
\newblock {On stable fields and weight}.
\newblock {\em J.\ Inst.\ Math.\ Jussieu 10}, 2 (2011), 349--358.

\bibitem{macintyre}
{\sc Macintyre, A.}
\newblock {On $\omega_1$-categorical theories of fields}.
\newblock {\em Fund.\ Math. 71}, 1 (1971), 1--25.

\bibitem{marker}
{\sc Marker, D.}
\newblock {\em {Model Theory: An Introduction}}.
\newblock {Springer}, 2002.
\newblock {Graduate Texts in Mathematics 217}.

\bibitem{markervf}
{\sc Marker, D.}
\newblock {Model Theory of Valued Fields}.
\newblock {\textit{Lecture notes, Available at}
  \href{http://homepages.math.uic.edu/~marker/valued_fields.pdf}{\ttfamily
  homepages.math.uic.edu/\\$\sim$marker/valued\_fields.pdf}}, 2018.

\bibitem{presteltop}
{\sc Prestel, A.}
\newblock {Algebraic number fields elementarily determined by their absolute
  Galois group}.
\newblock {\em {Israel J.\ Math.} 73}, 2 (1991), 199--205.

\bibitem{prestelandroq}
{\sc Prestel, A., and Roquette, P.}
\newblock {\em {Formally $p$-adic Fields}}.
\newblock {Springer}, 1984.
\newblock {Lecture Notes in Mathematics 1050}.

\bibitem{prestelziegler}
{\sc Prestel, A., and Ziegler, M.}
\newblock {Model-theoretic methods in the theory of topological fields}.
\newblock {\em {J.\ Reine Angew.\ Math.} 299}, 300 (1978), 318--341.

\bibitem{undecrob}
{\sc Robinson, R.}
\newblock Undecidable rings.
\newblock {\em Trans.\ Amer.\ Math.\ Soc. 70}, 1 (1951), 137--159.

\bibitem{shelahndep}
{\sc Shelah, S.}
\newblock {Strongly dependent theories (Sh:863)}.
\newblock {\em Israel J. Math. 204\/} (2014), 1--83.

\bibitem{shlapentokhvidela}
{\sc Shlapentokh, A., and Videla, C.}
\newblock {Definability and decidability in infinite algebraic extensions}.
\newblock {\em Ann.\ Pure Appl.\ Logic 165}, 7--8 (2014), 1243--1262.

\bibitem{shoenfield}
{\sc Shoenfield, J.~R.}
\newblock {\em {Recursion Theory}}.
\newblock Springer-Verlag, 1993.
\newblock {Lecture Notes in Logic, Volume 1}.

\bibitem{simon}
{\sc Simon, P.}
\newblock {\em {A Guide to NIP Theories}}.
\newblock Cambridge University Press, 2015.
\newblock {Lecture Notes in Logic}.

\bibitem{tmr}
{\sc Tarski, A., Mostowski, A., and Robinson, R.}
\newblock {\em {Undecidable Theories}}.
\newblock North-Holland, 1953.
\newblock {Studies in Logic and the Foundations of Mathematics}.

\bibitem{btyrrelpac}
{\sc Tyrrell, B.}
\newblock {Finite Undecidability in Fields II: PAC, PRC \& P$p$C Fields}.
\newblock {\em \emph{\href{https://arxiv.org/abs/2212.12918}{\ttfamily
  arXiv:2212.12918}}\/} (2023).

\bibitem{tyrrellphd}
{\sc Tyrrell, B.}
\newblock {\em {\it Undecidability in some Field Theories}}.
\newblock PhD thesis, {University of Oxford}, 2023.

\bibitem{ziegler}
{\sc Ziegler, M.}
\newblock {Einige unentscheidbare K\"{o}rpertheorien}.
\newblock {\em Enseign.\ Math.\ II 28}, 1--2 (1982), 269--280.

\end{thebibliography}
\bibliographystyle{acm-new} 
\end{document}